\documentclass[a4paper,11pt]{amsart}

\usepackage[hang,flushmargin]{footmisc}

\usepackage{orcidlink}

\usepackage{graphicx, rotating}
\usepackage{subfigure}
\usepackage[font=footnotesize,labelfont=footnotesize]{caption}
\usepackage[utf8]{inputenc}
\usepackage[toc]{appendix}
\usepackage{twoopt}

\usepackage{amsthm}
\usepackage{amsmath,a4wide}

\usepackage{tikz}
\usetikzlibrary{shapes.geometric}

\makeatletter
\newcommand*{\trinum}{}
\DeclareRobustCommand*{\trinum}[1]{%
  \ensuremath{%
    \mathpalette\@trinum{#1}%
  }%
}
\newdimen\trinum@sep
\newdimen\trinum@rule
\newcommand*{\@trinum}[2]{%
  \settowidth\trinum@sep{$\m@th#1\mkern1mu$}%
  \setlength{\trinum@rule}{.8\trinum@sep}
  \tikz\node[
    regular polygon,
    regular polygon sides=3,
    draw,
    line width=\trinum@rule,
    inner sep=\trinum@sep,
    scale=.5,
  ]{$\m@th#1#2$};
}
\makeatother

\usepackage{amssymb}
\numberwithin{equation}{section}
\usepackage{hyperref}
\usepackage{mathtools}
\usepackage{amsfonts}
\usepackage[foot]{amsaddr}
\usepackage{commath}

\usepackage{marvosym}
\usepackage{verbatim}

\usepackage{bookmark}
\bookmarksetup{depth=2}

\usepackage{nameref}

\usepackage{enumerate}
\usepackage{pdfsync}
\usepackage{caption}
\usepackage{hyperref}
\usepackage{enumerate}
\mathtoolsset{showonlyrefs}

\theoremstyle{plain}
\newtheorem{theorem}{Theorem}[section]

\theoremstyle{plain}
\newtheorem*{theorem*}{Theorem}

\theoremstyle{plain}

\theoremstyle{plain}

\newtheorem*{lemma*}{Lemma}
\theoremstyle{plain}

\newtheorem{proposition}[theorem]{Proposition}
\theoremstyle{definition}

\theoremstyle{remark}

\theoremstyle{remark}

\theoremstyle{definition}

\theoremstyle{plain}
\newtheorem{conjecture}[theorem]{Conjecture}
\theoremstyle{definition}

\providecommand{\norm}[1]{\left\lVert #1 \right\rVert}

\newcommand{\R}{\mathbb{R}}
\newcommand{\C}{\mathbb{C}}

\newcommand{\Rd}{\mathbb{R}^d}
\newcommand{\Z}{\mathbb{Z}}
\newcommand{\N}{\mathbb{N}}
\newcommand{\Lt}[1][d]{L^2(\R^{#1})}
\newcommand{\G}{\mathcal{G}}
\newcommand{\F}{\mathcal{F}}

\renewcommand{\l}{\lambda}
\renewcommand{\L}{\Lambda}
\newcommand{\LHex}{\Lambda_2} 

\newcommand{\vol}{\textnormal{vol}}

\newcommandtwoopt{\xarrow}[2][0.5cm][0]{\mathrel{\rotatebox[origin=c]{#2}{$\xrightarrow{\rule{#1}{0pt}}$}}}

\newcommand{\agm}{\mathsf{ag}_2}
\newcommand{\agc}{\mathsf{ag}_3}

\begin{document}

\title[AGM of Gauss, Ramanujan's corresponding theory, and spectral bounds]{
		The AGM of Gauss, Ramanujan's corresponding theory, and spectral bounds of self-adjoint operators
}

\author[M.\ Faulhuber]{Markus Faulhuber \orcidlink{0000-0002-7576-5724}
}
\author[A.\ Gumber]{Anupam Gumber \orcidlink{0000-0001-5146-3134}
}
\author[I.\ Shafkulovska]{Irina Shafkulovska \orcidlink{0000-0003-1675-3122}
}
\address{NuHAG, Faculty of Mathematics, University of Vienna,\newline Oskar-Morgenstern-Platz 1, 1090 Vienna, Austria}
\email{markus.faulhuber@univie.ac.at}
\email{anupam.gumber@univie.ac.at}
\email{irina.shafkulovska@univie.ac.at}

\thanks{The authors were supported by the Austrian Science Fund (FWF) projects P33217 and TAI6.}

\subjclass[2020]{33C05, 33C67, 33C80, 42C15}
\keywords{arithmetic-geometric mean, Gabor system, theta function, lattice, spectral bounds}

\begin{abstract}
	We study the spectral bounds of self-adjoint operators on the Hilbert space of square-integrable functions, arising from the representation theory of the Heisenberg group. Interestingly, starting either with the von Neumann lattice or the hexagonal lattice of density 2, the spectral bounds obey well-known arithmetic-geometric mean iterations. This follows from connections to Jacobi theta functions and Ramanujan's corresponding theories. As a consequence we re-discover that these operators resemble the identity operator as the density of the lattice grows. We also prove that the conjectural value of Landau's constant is obtained as the cubic arithmetic-geometric mean of $\sqrt[3]{2}$ and 1, which we believe to be a new result.
\end{abstract}

\vspace*{-1cm}
\maketitle


%
%
%
%
%

\section{Introduction}\label{sec:Intro}
\noindent
In this work we connect the problem of estimating the spectral bounds of certain self-adjoint operators on $\Lt[]$ with the arithmetic-geometric mean of Gauss and Ramanujan's corresponding theory. The operators involved arise from Gaussian Gabor systems over scaled von Neumann lattices and hexagonal lattices. To some extent, we can also draw conclusions for rectangular lattices. A common theme for all these lattices is that they contain certain root systems and, in fact, we cover all possible root lattices in dimension 2. Gauss has already shown that the squares of the classical Jacobi theta functions, which are used to define the elliptic modulus of complete elliptic integrals of the first kind, obey the arithmetic-geometric mean process. Moreover, it has been described by Mumford \cite{Mum_Tata_I} how theta functions play an important role in the representation theory of the Heisenberg group. Interestingly, the specific spectral problem also arises from unitary representations of the Heisenberg group and is of importance in quantum mechanics and communication theory. To the best of our knowledge, a connection between the arithmetic-geometric mean iteration and the spectral problem has so far neither been observed nor been studied in the literature.

Restricting (the squares of) Jacobi's theta functions to certain arguments allows us to connect to the spectral bounds for Gaussian Gabor frame operators over scaled von Neumann lattices, i.e., square lattices. We can then use a theory developed by Gauss to derive our first result. Moreover, Ramanujan's corresponding theory for elliptic functions enters the scene and provides a cubic analogue to the arithmetic-geometric mean. This has consequences for our considered spectral problem over the hexagonal lattice and gives our second result.

We denote the Gaussian function of $\Lt[]$ unit norm by $\varphi(t) = 2^{1/4} e^{-\pi t^2}$. For an element $(x,\omega, \tau)$ from the (polarized) Heisenberg group $\mathbf{H}$ (see \cite[Chap.\ 1]{Fol89} or \cite[Chap.\ 9]{Gro01}), we denote the unitary operator arising from its Schrödinger representation by $\pi(x,\omega;\tau)$. We will only be interested in $\pi(x,\omega;0)$, which we will simply denote by $\pi(x,\omega)$.

\noindent
A Gaussian Gabor system over a lattice $\L \subset \R^2$ is a structured function system of the form
\begin{equation}
	\G(\varphi,\L) = \{ \pi(\l) \varphi \mid \l \in \L\}, \quad \l = (x,\omega) \in \R^2.
\end{equation}
The associated self-adjoint Gabor frame operator acts on functions $f \in \Lt[]$ by the rule
\begin{equation}
	S_\L f = \vol(\L) \sum_{\l \in \L} \langle f, \pi(\l) \varphi \rangle \, \pi(\l) \varphi.
\end{equation}
Our main results, Theorem \ref{thm:agm2} and Theorem \ref{thm:agm3}, concern the sharp spectral bounds of $S_\L$. The spectral bounds of the operator will be denoted by $A_\L$ and $B_\L$ (we will later also pass the density of the lattice $\L$ as argument to the bounds), hence,
\begin{equation}
	A_\L \norm{f}_2^2 \leq \langle S_\L f, f \rangle \leq B_\L \norm{f}_2^2, \quad \forall f \in \Lt[].
\end{equation}
We denote the scaled von Neumann lattice and hexagonal lattice of density $\alpha > 0$ by
\begin{equation}\label{eq:lattices}
	\L_{1\times1}(\alpha) = \alpha^{-1/2} \Z \times \alpha^{-1/2} \Z
	\quad \text{and} \quad
	\L_2(\alpha) = \alpha^{-1/2} \sqrt{2 / \sqrt{3}}
	\begin{pmatrix}
		1 & \frac{1}{2}\\
		0 & \frac{\sqrt{3}}{2}
	\end{pmatrix}
	\Z^2,
	\quad \text{respectively}.
\end{equation}

\subsection{Main results}
\begin{theorem}[AGM2]\label{thm:agm2}
	Let $\L_{1\times1}(2^{n}) = 2^{-n/2} \Z \times 2^{-n/2} \Z$, $n \in \N$, denote the scaled von Neumann lattice of density $2^n$. Denote the spectral bounds of the frame operator by $A_{1\times1}(2^n)$ and $B_{1\times1}(2^n)$. Then, they obey the arithmetic-geometric mean iteration:
	\begin{equation}
		B_{1\times1}(2^{n+1}) = \frac{A_{1\times1}(2^n)+B_{1\times1}(2^n)}{2}
		\quad \text{ and } \quad
		A_{1\times1}(2^{n+1}) = \sqrt{A_{1\times1}(2^n) B_{1\times1}(2^n)}.
	\end{equation}
	The constants $A_{1\times1}(2)$ and $B_{1\times1}(2)$ satisfy the relation $B_{1\times1}(2)/A_{1\times1}(2) = \sqrt{2}$. Furthermore, denoting the classical arithmetic-geometric mean by $\agm$, we have
	\begin{equation}
		A_{1\times1}(2) = G = \frac{1}{\agm (\sqrt{2}, 1)} = \theta_4(e^{-\pi})^2 = 0.8346 \ldots,
	\end{equation}
	where $G$ denotes Gauss' constant and $\theta_4$ a Jacobi theta constant.
\end{theorem}

\begin{theorem}[AGM3]\label{thm:agm3}
	Let $\LHex(2\cdot 3^{n-1}) = 2^{-1/2}3^{-(n-1)/2} \LHex$, $n \in \N$, denote the scaled hexagonal lattice of density $2\cdot 3^{n-1}$ and the respective spectral bounds of the frame operator by $A_2(2\cdot 3^{n-1})$ and $B_2(2\cdot 3^{n-1})$. Then, they obey the cubic arithmetic-geometric mean iteration: 
	\begin{equation}
		B_2(2\cdot 3^{n}) = \frac{B_2(2\cdot 3^{n-1}) +2A_2(2\cdot 3^{n-1})}{3}
	\end{equation}
 \centerline{ and }
 	\begin{equation}
		A_2(2\cdot 3^{n})  = \sqrt[3]{A_2(2\cdot 3^{n-1}) \, \frac{A_2(2\cdot 3^{n-1}) ^2+A_2(2\cdot 3^{n-1}) B_2(2\cdot 3^{n-1}) +B_2(2\cdot 3^{n-1}) ^2}{3}} \, .
	\end{equation}
	The constants $A_2(2)$ and $B_2(2)$ satisfy the relation $B_2(2)/A_2(2)  = \sqrt[3]{2}$. Furthermore, denoting the cubic arithmetic-geometric mean by $\agc$, we have
	\begin{equation}
		A_2(2) = \frac{1}{2 \mathcal{L}_+} = \frac{1}{\agc(\sqrt[3]{2}, 1)} = b(e^{-\frac{2 \pi}{\sqrt{3}}}) = 0.920371 \ldots,
	\end{equation}
	where $\mathcal{L}_+$ is the conjectural value of Landau's constant and $b$ the cubic pendant to $\theta_4^2$.
\end{theorem}

\noindent
The constant $A_2(2)$ may be referred to as an equianharmonic constant \cite{AbrSte72}, \cite{Fau21-RJ}, \cite[Chap.\ 7]{Fin03}.

\subsection*{Remark}\label{sec:remark}
The connection of the spectral bounds of the Gaussian frame operator to Gauss' constant and the conjectural value of Landau's constant is somewhat mysterious and simply pops out of the computations. Also, to the best of our knowledge the equality
\begin{equation}\label{eq:Landau_ag3}
	\agc(\sqrt[3]{2}, 1) = 2 \mathcal{L}_+,
\end{equation}
seemed to be unknown so far and we will provide a short proof for this equation as well.

\section{Root systems and lattices}\label{sec:root}
\begin{figure}[htp!]
	\subfigure[$\mathsf{A}_1 \times \mathsf{A}_1$]{
		\includegraphics[width=.225\textwidth]{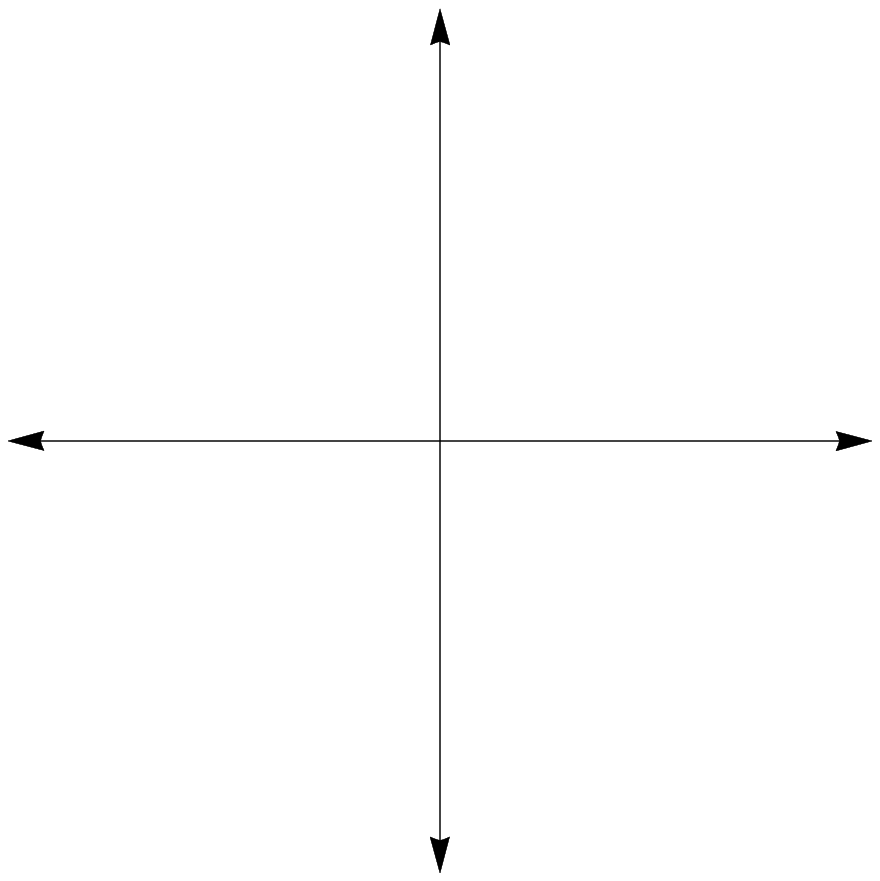}
	}
	\hfill
	\subfigure[$\mathsf{B}_2$]{
		\includegraphics[width=.225\textwidth]{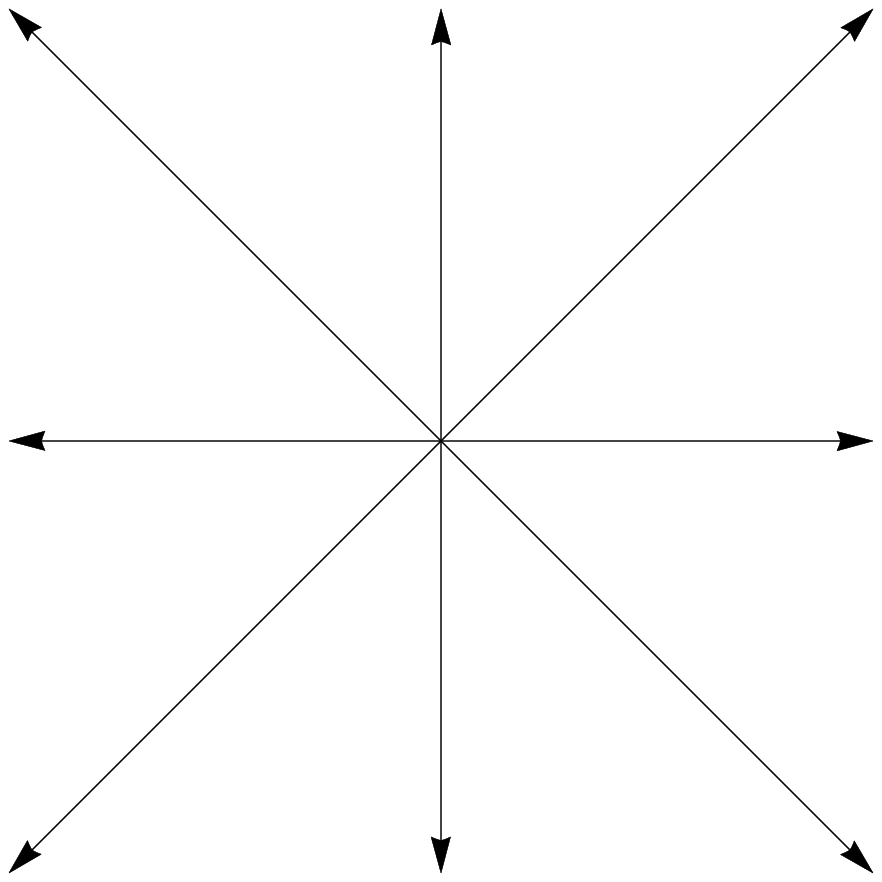}
	}
	\hfill
	\subfigure[$\mathsf{A}_2$]{
		\includegraphics[width=.225\textwidth]{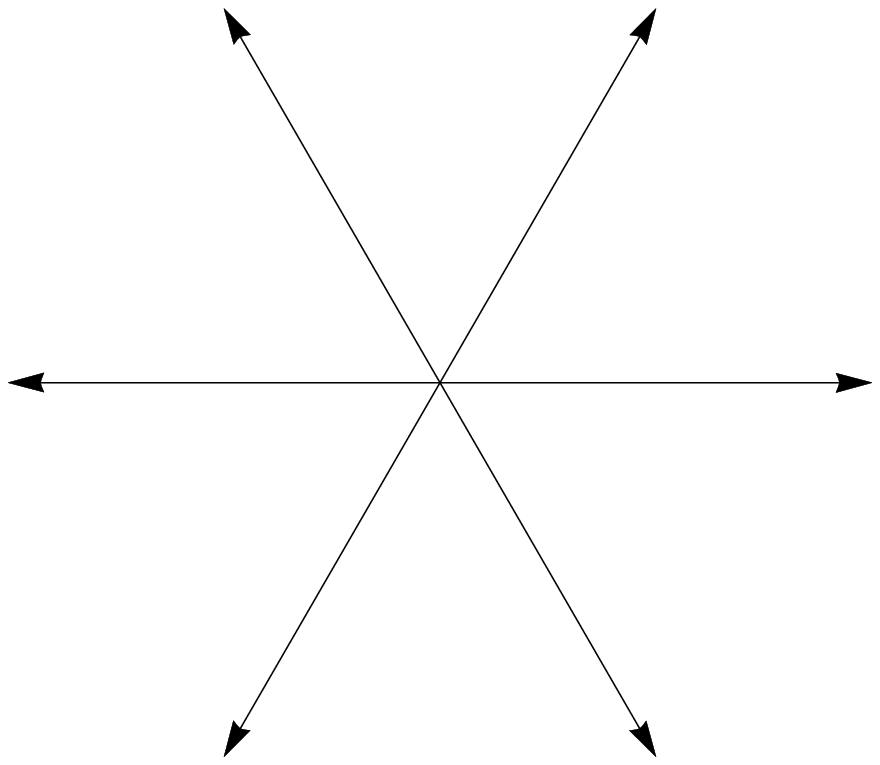}
	} 
	\\
	\subfigure[$\mathsf{D}2$]{
		\includegraphics[width=.225\textwidth]{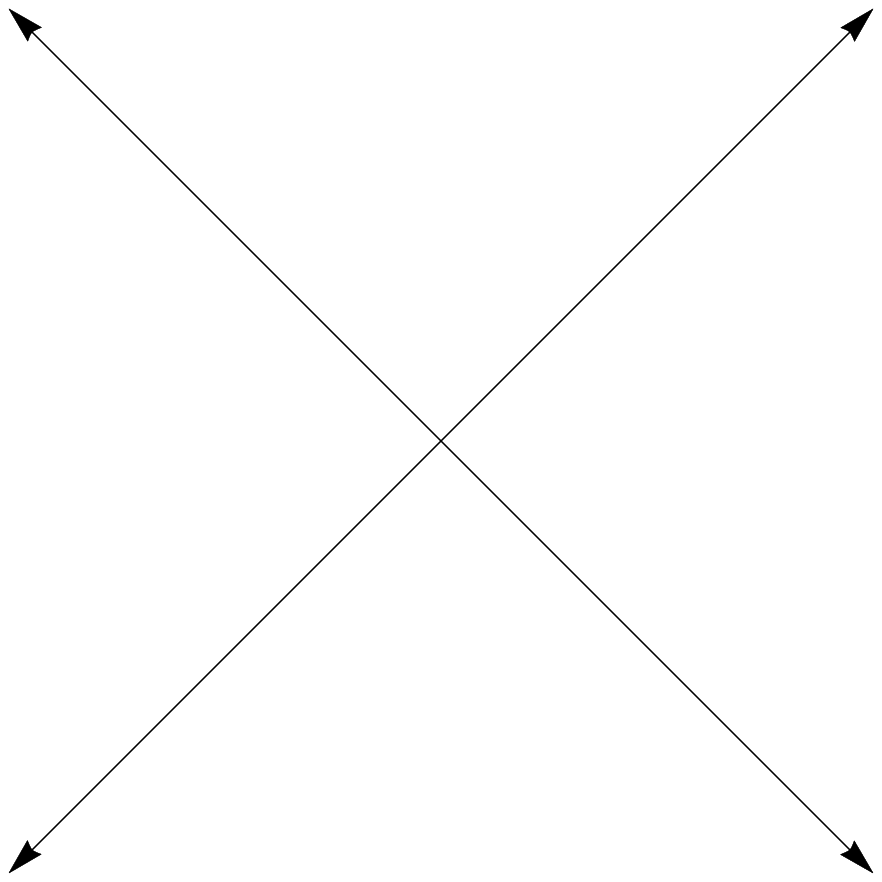}
	}
	\hfill
	\subfigure[$\mathsf{C}_2$]{
		\includegraphics[width=.225\textwidth]{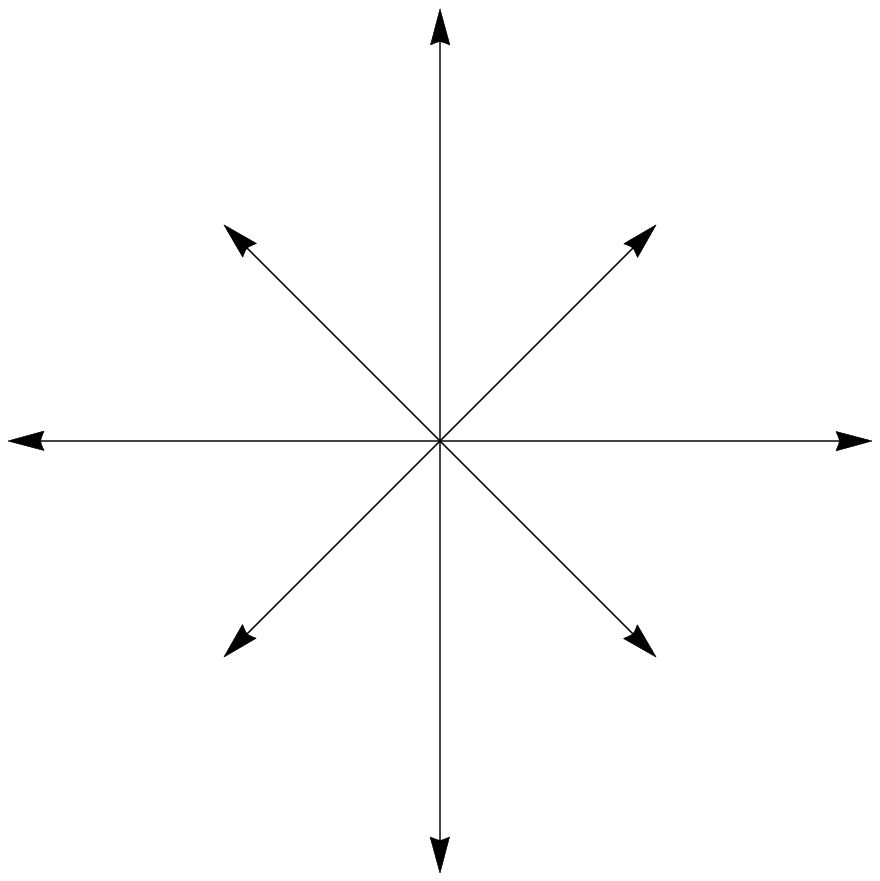}
	}
	\hfill
	\subfigure[$\mathsf{G}_2$]{
		\includegraphics[width=.225\textwidth]{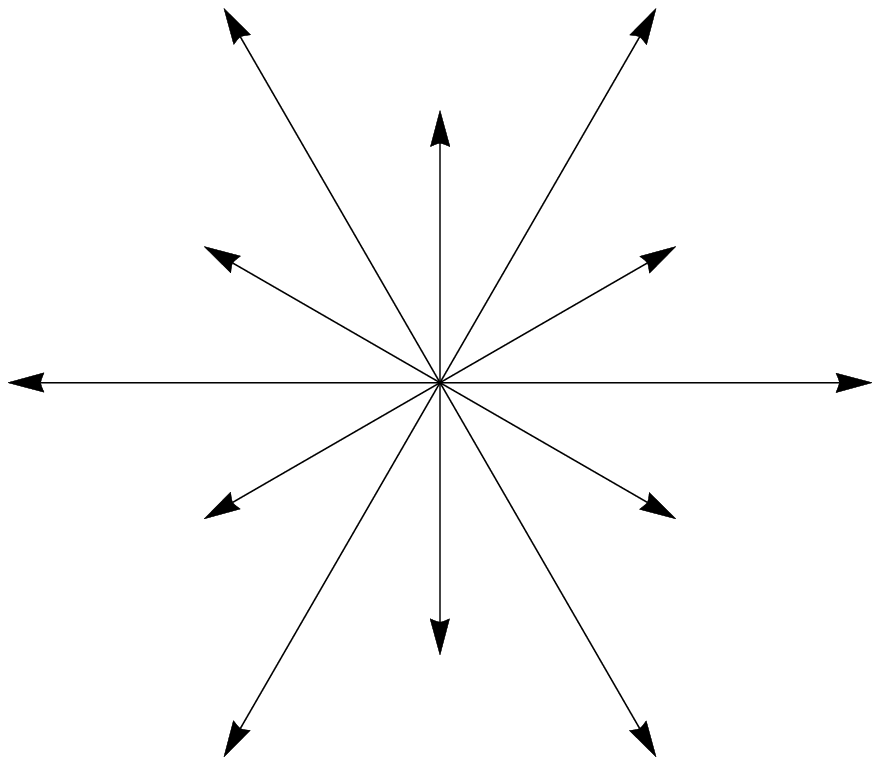}
	}
	\caption{The root systems $\mathsf{A}_1 \times \mathsf{A}_1$ and $\mathsf{D}_2$ are isomorphic as well as the root systems $\mathsf{B}_2$ and $\mathsf{C}_2$. All of them generate a (scaled) von Neumann lattice, by considering all integer linear combinations. The other existing root systems are $\mathsf{A}_2$ and $\mathsf{G}_2$. Both generate a hexagonal lattice. Note that $a A_1 \times b A_1$, $a, b > 0$ is also a root system which gives a rectangular lattice.}\label{fig:root_systems}
\end{figure}

\noindent
A root system is a finite set of vectors in Euclidean space $\Rd$ with exceptionally high symmetries. As we are only concerned with $\R^2$ in this work, we start with a full list of all possible root systems in $\R^2$ together with their illustrations in Figure \ref{fig:root_systems}. 

\vspace*{-0.2cm}
\noindent
By $v_1 \cdot v_2$ we denote the Euclidean inner product of the two vectors. A finite set $R$ of vectors is called a root system if
\begin{enumerate}[(i)]
	\item $0 \notin R$ and $\text{span}(R) = \Rd$.
	\item If $v_1 \in R$, called a root, then $-v_1 \in R$ and if for $r \in \R$ we have $r \, v_1 \in R$, then $r=\pm1$.
	\item If $v_1, v_2 \in R$, then
	\begin{equation}
		v_2 - 2 \frac{v_1 \cdot v_2}{v_1 \cdot v_1} \, v_1 \in R.
	\end{equation}
	\item If $v_1, v_2 \in R$, then
	\begin{equation}
		2 \frac{v_1 \cdot v_2}{v_1 \cdot v_1} \in \Z.
	\end{equation}
\end{enumerate}
We note that in dimension 1, there is only one root system, which is called $\mathsf{A}_1$ and only contains $\{\pm1\}$ (up to scaling). It is contained in the only 1-dimensional lattice, which is  $\Z$.

\begin{figure}[h!tp]
	\subfigure[The von Neumann lattice $\Z^2$ and the root system $\mathsf{A}_1 \times \mathsf{A}_1$]{
		\includegraphics[width=.4\textwidth]{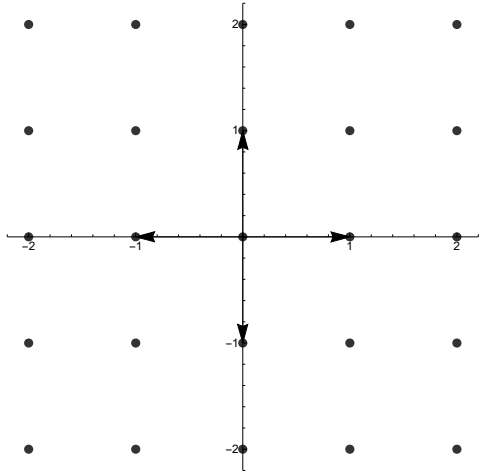}
	}
	\hfill
	\subfigure[The scaled and rotated von Neumann lattice of density 2 and the root system $\mathsf{C}_2$]{
		\includegraphics[width=.4\textwidth]{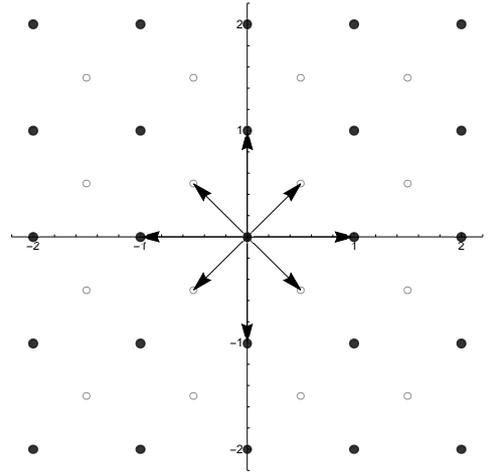}
	}
	\caption{The von Neumann lattice $\Z^2$ contains the root system $\mathsf{A}_1 \times \mathsf{A}_1$. By adding new points in the center of the fundamental cell (deep hole) we obtain a lattice with twice the density. The original lattice is contained as a sub-lattice and the new lattice contains the root system $\mathsf{C}_2$. The new lattice is merely a rotation of the scaled original by 45 degrees.}\label{fig:lattices_vN}
\end{figure}

\begin{figure}[h!tp]
	\subfigure[The hexagonal lattice of density 1 and the root system $\mathsf{A}_2$]{
		\includegraphics[width=.4\textwidth]{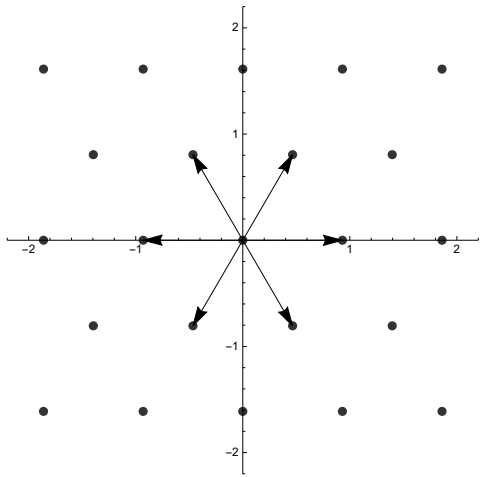}
	}
	\hfill
	\subfigure[The hexagonal lattice of density 3 and the root system $\mathsf{G}_2$]{
		\includegraphics[width=.4\textwidth]{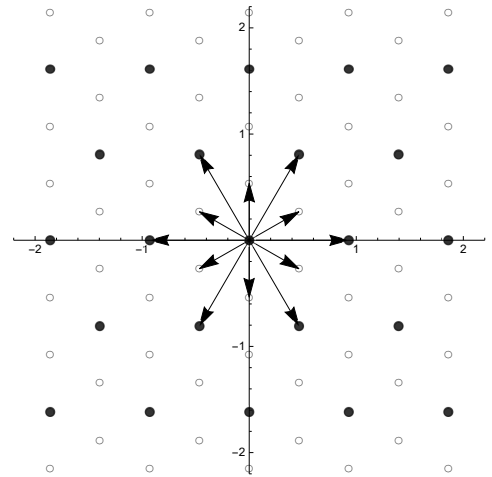}
	}
	\caption{The hexagonal lattice of density 1 contains the root system $\mathsf{A}_2$. By adding new points in the center of the fundamental triangle (deep hole) we obtain a lattice with thrice the density. The original lattice is contained as a sub-lattice and the new lattice contains the root system $\mathsf{G}_2$. The new lattice is merely a rotation of the scaled original by 30 degrees.}\label{fig:lattices_A2}
\end{figure}

The 6 root systems in dimension 2 can be used to construct lattices, or, phrased somewhat differently, are contained in certain lattices. The only lattices which contain root systems are the (scaled) von Neumann lattice, rectangular lattices and the hexagonal lattice. This is illustrated in Figure \ref{fig:lattices_vN} and Figure \ref{fig:lattices_A2}.

In this work, we are only interested in lattices which contain a root system, which reduces our interest mainly to the von Neumann lattice and the hexagonal lattice, and to some extent to rectangular lattices. We will use the following notation for these lattice, which is closely related to the $\mathsf{A}_*$ notation for root systems:
\begin{itemize}
	\item scaled von Neumann lattice: $\L_{1\times1}(\alpha) = \alpha^{-1/2} \Z^2$,
	\item hexagonal lattice: $\L_2(\alpha) = \alpha^{-1/2} \sqrt{\frac{2}{\sqrt{3}}}
	\begin{pmatrix}
		1 & \frac{1}{2}\\
		0 & \frac{\sqrt{3}}{2}
	\end{pmatrix} \Z^2$,
	\item rectangular lattice: $\L_{a \times b}(\alpha) = \alpha^{-1/2} \left(a \Z \times b \Z\right)$, with $a b = 1$.
\end{itemize}
In the above notation, we made the density $\alpha$ of the lattice explicit. The density is the average number of lattice points per unit area (more details are given in Section \ref{sec:Gabor} below). Sometimes, we will also suppress the notation and simply write $\L_{1\times1}$, $\L_2$ and $\L_{a \times b}$, for the scaled von Neumann, hexagonal and rectangular lattices, respectively. The density, however, plays an important role in our main results, which is why we then prefer the above notation.

\section{Gabor systems over lattices}\label{sec:Gabor}
\noindent
We consider the Hilbert space of square-integrable functions on the line, denoted by $\Lt[]$. The inner product and the norm are given by
\begin{equation}
	\langle f, g \rangle = \int_\R f(t) \overline{g(t)} \, dt
	\quad \text{ and } \quad
	\norm{f}_2^2 = \langle f, f \rangle,
	\quad \text{ respectively}.
\end{equation}
The unitary operators of translation (time-shift operator) and modulation (frequency-shift operator) are given by
\begin{equation}
	T_x f(t) = f(t-x)
	\quad \text{ and } \quad
	M_\omega f(t) = e^{2 \pi i \omega t},
	\quad \text{ satisfying } \quad
	M_\omega T_x =  e^{2 \pi i \omega x} T_x M_\omega,
\end{equation}
respectively. We recall that $\varphi(t) = 2^{1/4} e^{-\pi t^2}$ is the normalized Gaussian function. The Fourier transform of a (suitably nice) function on the line is given by
\begin{equation}
	\F f(\omega) = \int_\R f(t) e^{-2 \pi i \omega t} \, dt.
\end{equation}
Note that the Gaussian is invariant under the Fourier transform, i.e., it is an eigenfunction with eigenvalue 1: $\F \varphi = \varphi$. This is also true for the multi-variate Gaussians $\varphi \otimes \varphi$ and the Fourier transform over $\R^2$. We mention this as we will frequently use the Poisson summation formula (introduced in Section \ref{sec:lattice_theta}) with (versions of) a Gaussian of the form $\varphi \otimes \varphi$.

Combining the action of translation and modulation gives a so-called time-frequency shift, which also appears as the Schrödinger representation of an element of the (polarized) Heisenberg group (\cite[Chap.\ 1]{Fol89}, \cite[Chap.\ 9]{Gro01});
\begin{equation}
	\pi(z) = \pi(x,\omega), \quad z=(x,\omega) \in \R^2.
\end{equation}
In this context we refer to $\R^2$ as the time-frequency plane. As the (polarized) Heisenberg group $\mathbf{H}$ is a non-Abelian group, it follows that time-frequency shifts do not commute in general. This can also be seen from the fact that already $T_x$ and $M_\omega$ do not commute in general. For a nice treatise on the role of the Heisenberg group in harmonic analysis we refer to \cite{How80}. We will let time-frequency shifts act on the Gaussian function $\varphi$:
\begin{equation}
	\pi(z) \varphi(t) = M_\omega T_x \varphi(t) = \varphi(t-x) e^{2 \pi i \omega t} = 2^{1/4} e^{-\pi (x-t)^2} e^{2 \pi i \omega t}, \quad t \in \R.
\end{equation}
A question, going back to von Neumann in quantum mechanics \cite{Neumann_Quantenmechanik_1932}, and studied later independently by Gabor for the purpose of communication theory \cite{Gab46}, is whether the set
\begin{equation}
	\G_0(\varphi, \Z^2) = \{ \pi(k,l) \varphi \mid (k,l) \in \Z^2\}
\end{equation}
is complete in $\Lt[]$. Therefore, the lattice $\Z^2$ is often referred to as the von Neumann lattice in the context of quantum mechanics. As we now know, this is indeed the case \cite{BarButGirKla71}, \cite{Per71} (see also \cite{GroHaiRom16}). We also refer to \cite{BooZakZuc83} for a treatise on the (over)completeness of coherent states over the von Neumann lattice and its connections to zeros of theta functions.

In time-frequency analysis, one is interested in stable expansions of the form
\begin{equation}\label{eq:expansion}
	f = \sum_{(k,l) \in \Z^2} c_{k,l} \, \pi(k,l) \varphi,
\end{equation}
where $(c_{k,l}) \in \ell^2(\Z^2)$. For the above system $\G_0(\varphi, \Z^2)$ this is not possible, due to a manifestation of the uncertainty principle: the Balian-Low theorem \cite{Bal81}, \cite{Low85}. We will elaborate on this fact a bit later, at the end of this section.

In the above expansion \eqref{eq:expansion}, the Gaussian function $\varphi$ may be replaced by any other (suitably nice) function $g \in \Lt[]$ and the index set $\Z^2$ by a lattice $\L \subset \R^2$. A lattice in the time-frequency plane is a discrete co-compact subgroup of $\R^2$, i.e., the integer span of a basis for $\R^2$. We can write any lattice as
\begin{equation}
	\L = M \Z^2 = (v_1,v_2)\Z^2 = \{k v_1 + l v_2 \mid (k,l) \in \Z^2 \}.
\end{equation}
The matrix $M$ contains the column vectors $v_1$ and $v_2$ which constitute a basis for $\R^2$, in other words, $M \in \mathrm{GL}(2,\R)$. We write $\vol(\L)$ for the co-volume of the lattice $\L$, which is given by
\begin{equation}
	\vol(\L) = |\det(M)|,
	\quad \text{ for } \quad
	\L = M \Z^2.
\end{equation}
This leads to the study of Gabor systems, also known as Weyl-Heisenberg systems in mathematical physics \cite{Gos11}, of the form:
\begin{equation}
	\G(g, \L) = \{ \pi(\l) g \mid \l \in \L \}.
\end{equation}
The question for stable expansions of a function with respect to the Gabor system is equivalent to determining whether $\G(g, \L)$ is a Gabor frame for $\Lt[]$. The system $\G(g,\L)$ is a frame if and only if there exist positive constants $0<A \leq B < \infty$ (depending on $g$ and $\L$) such that
\begin{equation}\label{eq:frame}
	A \norm{f}_2^2 \leq \vol(\L) \, \sum_{\l \in \L} |\langle f, \pi(\l) g \rangle|^2 \leq B \norm{f}_2^2, \quad \forall f \in \Lt[].
\end{equation}
This is just a spectral inequality for the self-adjoint Gabor frame operator, denoted by $S_{g,\L}$:
\begin{equation}\label{eq:frame_expansion}
	S_{g,\L} f = \vol(\L) \, \sum_{\l \in \L} \langle f, \pi(\l) g \rangle \, \pi(\l) g,
	\quad \text{ so, } \quad
	A \norm{f}_2^2 \leq \langle S_{g,\L} f, f \rangle \leq B \norm{f}_2^2, \quad \forall f \in \Lt[].
\end{equation}
Note that we have slightly adjusted the standard formulation of the frame operator and the frame inequality to an equivalent statement by multiplying by the co-volume of the lattice. This will turn out to be a convenient normalization for our main results. Also, note that the sharpest possible constants in \eqref{eq:frame} are actually the spectral bounds of the operator $S_{g,\L}$:
\begin{equation}
	A = \norm{S_{g,\L}^{-1}}_{\Lt[] \to \Lt[]}^{-1}, \quad B = \norm{S_{g, \L}}_{\Lt[] \to \Lt[]}.
\end{equation}

We will now briefly come back to the above mentioned manifestation of the uncertainty principle in time-frequency analysis. The quantity $\vol(\L)^{-1}$ is called the density of the lattice and gives the average number of points (or available information) per unit area. For Gaussian Gabor systems of the form $\G(\varphi, \L)$ we simply write $S_\L$ for the frame operator. In this case we know that the necessary density condition $\vol(\L) < 1$, imposed by the Balian-Low theorem, is already sufficient for the Gabor system to produce a Gabor frame:
\begin{equation}
	\G(\varphi, \L) \text{ is a frame}
	\quad \Longleftrightarrow \quad
	\vol(\L) < 1.
\end{equation}
This is a consequence of the celebrated results of Lyubarskii \cite{Lyu92}, Seip \cite{Sei92_1}, Seip and Wallst\'en \cite{SeiWal92}. The condition $\vol(\L) < 1$ tells us that we need more than 1 time-frequency sample per unit area. This is comparable to the Nyquist rate in the classical Whittaker-Koltelnikov-Shannon sampling theorem for band-limited functions \cite{Kot33}, \cite{Sha49}, \cite{Whi15}.

\section{AGMs and lattice theta functions}\label{sec:agm}
\noindent
The arithmetic-geometric means (AGMs) of order $N$ are recursive constructions of sequences studied by Borwein and Borwein \cite{BorBor_Cubic_91}. The general iteration for $1 < N\in \N$ is initialized with two starting values $a_0$ and $b_0$, which are allowed to be complex numbers. We will only use the process for real values $a_0 > b_0 > 0$. Three sequences are defined recursively:
\begin{equation}
	a_{n+1}=\frac{a_n+(N-1)b_n}{N},
	\qquad c_{n+1}=\frac{a_n-b_n}{N},
	\qquad\text{where}\quad b_n^N=a_n^N-c_n^N.
\end{equation}
Then the sequences $\left( a_n\right)_{n\in\N}$ and $\left( b_n\right)_{n\in\N}$ converge to a common limit, denoted by $\mathsf{ag}_N(a,b)$, with convergence rate $N$, best justified by
\begin{equation}
a_{n+1}^N-b_{n+1}^N = \left(\frac{a_n-b_n}{N}\right)^N.
\end{equation}
We will consider the arithmetic-geometric means of order $2$ and $3$, with particular sequences. By substituting $c_n$, we can describe the $(a_n)_{n\in\N}$ and $(b_n)_{n\in\N}$ independently of $(c_n)_{n\in\N}$ as:
\begin{align}
	a_{n+1} &= \frac{a_n+b_n}{2} &\qquad\text{and}\qquad b_{n+1}&=\sqrt{a_n b_n}, &\qquad\qquad (N=2), \\
	a_{n+1} &= \frac{a_n+2b_n}{3} &\qquad\text{and}\qquad b_{n+1}&=\sqrt[3]{b_n\left(\frac{a_n^2+a_nb_n+b_n^2}{3}\right)}, &\qquad\qquad (N=3).
\end{align}
The case $N = 2$ simply leads to the classical arithmetic-geometric mean. This has already been studied intensively by Gauss, who proved that it can be used to exactly and efficiently compute elliptic integrals \cite[entry 98, entry 102]{Kle03} (see also \cite{Cox84} and \cite[pp.\ 446 ff.]{Gauss_WerkeIII}). Note that, for $x > 0$, it is obvious from the above iterations for $N=2$ and $N=3$ that $\agm$ and $\agc$ are homogeneous:
\begin{equation}
	\agm(x a_0, x b_0) = x \, \agm(a_0, b_0)
	\quad \text{ and } \quad
	\agc(x a_0, x b_0) = x \, \agc(a_0, b_0).
\end{equation}

\subsection{Jacobi theta functions}
Theta functions are classical objects, appearing in many branches of mathematics and also in other sciences. We refer to the textbooks of \cite{SteSha_Complex_03, WhiWat69} for proper introductions as well as to \cite[Chap.~4]{ConSlo98} for their connection with sphere packings and coverings and related topics. Still, we want to clarify the notation which we use in this work. We write the Jacobi theta functions in the following way;
\begin{align}
	\vartheta_2(z;q) = \sum_{k \in \Z} q^{(k+\frac{1}{2})^2} e^{2 \pi i (k+ \frac{1}{2}) z},
	\quad
	\vartheta_3(z;q) = \sum_{k \in \Z} q^{k^2} e^{2 \pi i k z},
	\quad
	\vartheta_4(z;q) = \sum_{k \in \Z} (-1)^k q^{k^2} e^{2 \pi i k z},
\end{align}
where $z \in \C$ and $0 < \abs{q} < 1$. It is also common to replace the nome $q$ by $e^{\pi i \tau}$. Then $\tau$ needs to be chosen from the Siegel upper half space
\begin{equation}
	\mathbb{H} = \{ z \in \C \mid \Im(z) > 0 \},
\end{equation}
in order to ensure convergence of the series. The $\vartheta$-functions are entire for $z \in \C$ and holomorphic for $\tau \in\mathbb{H}$. For our purposes, it suffices to consider the functions for $z = 0$, which carry the name theta-constants. We write
\begin{equation}
	\theta_m(q) = \vartheta_m(0,q), \quad m \in \{2,3,4\}.
\end{equation}
We remark that there is also the (prototype) Jacobi theta function
\begin{equation}
	\vartheta_1(z,q) = -i \sum_{k \in \Z} (-1)^k q^{(k+\frac{1}{2})^2} e^{2 \pi i (k + \frac{1}{2}) z}.
\end{equation}
It is an odd function in $z$: its theta constants vanishes identically and is not of interest here.

\subsection{The analogues}\label{sec:analogues}

In their paper on the cubic AGM \cite{BorBor_Cubic_91}, Borwein and Borwein introduced new $\theta$-like functions. Following the notation in \cite{BorBor_Cubic_91}, we set
\begin{align}
	a(q) & = \sum_{m,n \in \Z} q^{m^2+mn+n^2}\\
	b(q) & = \sum_{m,n \in \Z} \zeta^{n-m} q^{m^2+mn+n^2}
	\\
	c(q) & = \sum_{m,n \in \Z} q^{(m+\frac{1}{3})^2+(m+\frac{1}{3})(n+\frac{1}{3})+(n+\frac{1}{3})^2}
\end{align}
where $\zeta^3=1$, $\zeta\neq 1$ and $0 < \abs{q} < 1$. The analogy is understood through the correspondence 
\begin{equation}
	\begin{split}
		\theta_3(q)^2 \ \longleftrightarrow\ a(q), \qquad	\theta_4(q)^2 \  \longleftrightarrow \ b(q), \qquad
		\theta_2(q)^2 \ \longleftrightarrow \ c(q). 
	\end{split}
\end{equation}

\subsection{The AGM of theta functions}\label{sec:AGM_theta}
The triples of theta constants and their cubic analogues satisfy (see \cite[Chap.\ 4.4]{ConSlo98} for the classical theta functions and \cite{BorBor_Cubic_91} for the cubic analogues)
\begin{equation}\label{eq:correspondence}
\begin{aligned}
	2\theta_3(q^2)^2 &= \theta_3(q)^2+\theta_4(q)^2 &\longleftrightarrow&\qquad &3a(q^3) &= a(q) + 2b(q), && \\
	2\theta_2(q^2)^2& =\theta_3(q)^2 - \theta_4(q)^2 \qquad &\longleftrightarrow&\qquad  &3c(q^3) &= a(q) - b(q), && \\
	\theta_3(q)^4 &= \theta_4(q)^4+\theta_2(q)^4 \qquad&\longleftrightarrow&\qquad 	&a(q)^3&= b(q)^3 + c(q)^3.&& \\
\end{aligned}
\end{equation}
This shows that $\theta_3^2,\,\theta_4^2,\,\theta_2^2$ and $a,\,b,\,c$ fit within the $\agm$ and $\agc$ construction, respectively.
For more identities involving the cubic analogues we refer to the articles \cite{BorBorGar95} and \cite{HirGarBor94}.



\subsection{Lattice theta functions and Gaussian lattice sums}\label{sec:lattice_theta}
Recall that in order for the Gabor system $\G(\varphi, \L)$ to be a frame, the frame inequality has to be satisfied:
\begin{equation}\label{eq:frame_gaussian}
	A_\L \norm{f}_2^2 \leq \vol(\L) \, \sum_{\l \in \L} |\langle f, \pi(\l) \varphi \rangle|^2
	\leq B_\L \norm{f}_2^2, \quad \forall f \in \Lt[].
\end{equation}
Finding extremal functions $f$ (depending on $\L$) such that the sharp bounds are met is a cumbersome task. We may relax the problem by only considering the set $\{\pi(z) \varphi \mid z \in \R^2 \}$, which is known to be dense in $\Lt[]$ \cite[Chap.\ 1.5]{Gro01}. Noting that $\norm{\pi(z) \varphi}_2^2 = 1$, we have:
\begin{equation}\label{eq:frame_relaxed}
	A_\L \leq \vol(\L) \, \sum_{\l \in \L} |\langle \pi(z) \varphi, \pi(\l) \varphi \rangle|^2 \leq B_\L.
\end{equation}
A small computation shows that we have (see also \cite[Prop.\ (1.48)]{Fol89}, \cite[Chap.\ 1.5]{Gro01})
\begin{equation}
	|\left\langle \pi(z) \varphi, \pi(\l) \varphi\right\rangle|^2 = |\langle \varphi, \pi(\l-z) \varphi\rangle|^2	= e^{{-\pi}|\l-z|^2}.
\end{equation}
Since $\L$ has an additive group structure, we see that $\l \in \L$ if and only if $-\l \in \L$. So, we can actually re-write \eqref{eq:frame_relaxed} as
\begin{equation}
	A_\L \leq \vol(\L) \sum_{\l \in \L} e^{-\pi |\l+z|^2} \leq B_\L, \quad \forall z \in \R^2.
\end{equation}
It should be noted that the above inequality may actually not become sharp!

For a lattice $\L$ of unit density, i.e., $\vol(\L) = 1$, and $\alpha >0$, we now introduce the following family of lattice theta functions:
\begin{equation}
	\theta_{\L}(b;\alpha)= \sum_{\l \in \L} e^{-\pi\alpha|\l +b|^2}, \quad b \in \R^2.
\end{equation}
These lattice theta functions, which are Gaussian lattice sums shifted by $b=(b_1,b_2)$ in $\R^2$, will play a central role in this part. Actually, it is more their symplectic dual, the modulated Gaussian lattice sums, which will be of importance. For a lattice of unit density, this is
\begin{equation}
	\widehat{\theta}_\L (b;\alpha) = \sum_{\l \in \L} e^{-\pi \alpha |\l|^2} e^{2 \pi i \sigma(b, \l)}.
\end{equation}
At this point, it is necessary to introduce the symplectic form $\sigma(. \, , \, . )$ and the (symplectic) Poisson summation formula. For a suitable function $f$ the Poisson summation formula for a lattice $\L$ and its dual lattice $\L^\perp$ is (see \cite[Chap.\ 1]{Gro01})
\begin{equation}
	\sum_{\l \in \L} f(\l+x) = \vol(\L)^{-1} \sum_{\l^\perp \in \L^\perp} \F f(\l^\perp) e^{2 \pi i \l^\perp \cdot x}.
\end{equation}
We only give a characterization of the dual lattice in dimension 2, but the statement is easily transferred to higher dimensions. Denoting by $M^{-T}$ the inverse of $M$ transposed we have
\begin{equation}
	\L^\perp = \{ \l^\perp \in \R^2 \mid \l^\perp \cdot \l \in \Z, \, \forall \l \in \L \} = M^{-T} \Z^2, \quad \L = M \Z^2.
\end{equation}
As we are working in dimension 2, we can actually exploit the symplectic structure of the time-frequency plane (see \cite{Fol89}, \cite{Gos11}). We introduce the standard symplectic form $\sigma$, which is skew-symmetric and will replace the Euclidean inner produce in some computations;
\begin{equation}
	\sigma(z,z') = z_1 z_2' - z_2 z_1' = z \cdot \mathcal{J} z', \quad z=(z_1,z_2), \, z'=(z_1',z_2') \in \R^2.
\end{equation}
We denote by $\mathcal{J}$ the standard symplectic matrix. In $\R^2$ it is simply a rotation by 90 degrees;
\begin{equation}\label{eq:J}
	\mathcal{J} =
	\begin{pmatrix}
		0 & 1\\
		-1 & 0
	\end{pmatrix}.
\end{equation}
As we work in $\R^2$, we can use $\sigma$ to define the symplectic Fourier transform \cite{Gos11}, \cite{Gos17}:
\begin{equation}
	\F_\sigma F(z) = \iint_{\R^2} F(z') e^{2 \pi i \sigma(z,z')} \, dz'.
\end{equation}
The symplectic Fourier transform carries many properties of the ordinary Fourier transform. It is for example unitary. A main difference is that it is involutive, i.e.,
\begin{equation}
	\F_\sigma(\F_\sigma(F)) = F.
\end{equation}
Using the symplectic machinery, we can easily introduce a version of the Poisson summation formula for lattices in $\R^2$. We call it the symplectic Poisson summation formula:
\begin{equation}
	\sum_{\l \in \L} F(\l+z) = \vol(\L)^{-1} \sum_{\l^\circ \in \L^\circ} \F_\sigma F(\l^\circ) e^{2 \pi i \sigma(\l^\circ , z)}, \quad z \in \R^2.
\end{equation}
Here, $\L^\circ$ denotes the adjoint or symplectic dual lattice. This is just the usual dual lattice rotated by 90 degrees. It is characterized by commuting time-frequency shifts:
\begin{equation}\label{eq:adjoint_lattice}
	\L^\circ = \{\l^\circ \in \R^2 \mid \pi(\l)\pi(\l^\circ) = \pi(\l^\circ) \pi(\l), \; \forall \l \in \L \} = \mathcal{J} M^{-T} \Z^2.
\end{equation}
It should be noted that for 2-dimensional lattices the adjoint is actually simply a re-scaling of the lattice $\L$, i.e., for $\L(\alpha) = \alpha^{-1/2} M \Z^2$, with $M \in \mathrm{SL}(2,\R)$ it holds that
\begin{equation}
	\L^\circ = \alpha \L.
\end{equation}
We have the following functional equation, similar to the Jacobi identity, which follows from the symplectic Poisson summation formula:
\begin{equation}
	\theta_\L(b;\alpha) = \tfrac{1}{\alpha} \widehat{\theta}_\L(b;\tfrac{1}{\alpha}).
\end{equation}
The families of functions $\theta_\L(b;\alpha)$ and $\widehat{\theta}_\L(b;\alpha)$ have been studied thoroughly by Bétermin and Faulhuber for the special argument $\widetilde{b}=(1/2,1/2)$, where $b = M \, \widetilde{b}$, $\L = M \Z^2$, $M \in \mathrm{SL}(2,\R)$, is the center of the fundamental cell of the lattice \cite{BetFau21}. It should be evident that we use column vectors in $\R^2$, even though we write them as row vectors. Bétermin, Faulhuber and Steinerberger studied the case of $b$ being the minimizer of the lattice theta function \cite{BetFauSte21}. The subtlety in \cite{BetFauSte21} is that, in general, the minimizer depends on $\alpha$. In both cases, the hexagonal lattice turns out to be the global maximizer among lattices of unit density and all $\alpha > 0$, which gives a dual universal optimality result among lattices in the spirit of \cite{CohKum07}. The case $b=(0,0)$, which can be replaced by any lattice point and which is the maximizer of the lattice theta functions, has been fully treated by Montgomery \cite{Mon88}, proving that the hexagonal lattice is the unique minimizer in this case. In all cases, i.e., in \cite{BetFau21}, \cite{BetFauSte21}, \cite{Mon88}, the scaled von Neumann lattice is a critical point and, indeed, it is a saddle point in the set of lattices. These are few reasons why these lattices are of special interest. For more details on the connection to Gabor frames we refer to \cite{Faulhuber_Hexagonal_2018} and \cite{FauSte17}. For a more detailed discussion on the parametrization of lattices and symplectic methods for theta functions we refer to \cite{BetFau21}.

\subsection{The fundamental identity of Gabor analysis}
It is now advantageous to introduce the following notation, making the volume of the lattice explicit when needed:
\begin{equation}
	\L(\alpha) = \alpha^{1/2} M \Z^2, \quad M \in \mathrm{SL}(2,\R), \; \alpha > 0.
\end{equation}
We have the simple consequence that $\vol(\L(\alpha)) = \alpha.$ Note, that we can write $\L(\alpha) = \alpha^{1/2} \L(1)$. 
This provides the following significant number theoretic relation between Gaussian lattice sums and Gaussian Gabor frame bounds:
\begin{equation}\label{eq1:latticetheta}
	A_\L \leq \vol(\L(\alpha))\, \sum_{\l \in \L(\alpha)} e^{{-\pi}|\l +z|^2} = \alpha \, \theta_{\L(1)}(z;\alpha^{1/2}) \leq B_\L , \quad z \in \R^2.
\end{equation}

Next, we introduce the fundamental identity of Gabor analysis (FIGA), which is the Poisson summation formula in disguise (suppressing the dependency on the volume $\alpha$ again):
\begin{equation}\label{eq:FIGA}
	\vol(\L) \sum_{\l \in \L} \langle f_1, \pi(\l) g_1 \rangle \overline{\langle f_2, \pi(\l) g_2 \rangle} = \sum_{\l^\circ \in \L^\circ} \langle g_1, \pi(\l) g_2 \rangle \overline{\langle f_2, \pi(\l) f_1 \rangle}.
\end{equation}
We refer to \cite{FeiLue06}, \cite{GroKop19} or \cite{Jan98} for details and when the formula is applicable. As we use the Gaussian window $\varphi$, all requirements are, however, met. Using \eqref{eq:FIGA}, we get
\begin{align}
	\vol(\L) \, \sum_{\l \in \L} |\langle f, \pi(\l) \varphi \rangle|^2
	& = \vol(\L) \sum_{\l \in \L} \langle f, \pi(\l) \varphi \rangle \, \overline{\langle f, \pi(\l) \varphi \rangle}
	& & = \sum_{\l^{\circ} \in \L^{\circ}} \langle \varphi, \pi(\l^\circ) \varphi \rangle \, \overline{\langle f, \pi(\l^\circ) f \rangle}\\
	& \leq \sum_{\l^{\circ} \in \L^{\circ}} |\langle \varphi, \pi(\l^\circ) \varphi \rangle| \, |\overline{\langle f, \pi(\l^\circ) f \rangle}|
	& & \leq \sum_{\l^{\circ} \in \L^{\circ}} |\langle \varphi, \pi(\l^\circ) \varphi \rangle| \, \norm{f}_2^2.
\end{align} 
In the last step, we used the fact that $|\langle f, \pi(\l) f \rangle| \leq |\langle f, f \rangle|=\norm{f}_2^2$, which follows from the Cauchy-Schwarz inequality (see \cite[Lem.\ 4.2.1]{Gro01}). It readily follows (see also \cite{TolOrr95}) that  
\begin{equation}\label{eq:upperbound}
	B_\L \norm{f}_2^2 \leq \widetilde{B}_\L \norm{f}_2^2,
	\quad \text{ where } \quad
	\widetilde{B}_\L = \sum_{\l^{\circ} \in \L^{\circ}} |\langle \varphi, \pi(\l^\circ) \varphi \rangle| = \sum_{\l^\circ \in \L^\circ} e^{-\frac{\pi}{2} |\l^\circ|^2}.
\end{equation}
Hence, the quantity $\widetilde{B}_\L$ is a Bessel bound (not necessarily the sharpest one) for the Gabor system $\G(\varphi, \L)$ and we have $B_\L \leq \widetilde{B}_\L$. Combining observations \eqref{eq1:latticetheta} and \eqref{eq:upperbound}, we conclude:
\begin{equation}\label{eq2:latticetheta}
	A_\L \leq \vol(\L)\, \sum_{\l \in \L} e^{{-\pi}|\l +z|^2}= \alpha \, \theta_{\L(1)}(z;\alpha^{1/2}) \leq B_\L \leq \widetilde{B}_\L, \quad \forall z \in \R^2.
\end{equation}

\subsection{Ramanujan's Corresponding Theories}\label{sec:Ramanujan}
\noindent
In this part, we explain a nice connection between Gaussian lattice sums, lattice theta function and Ramanujan's corresponding theories of signature 2 and 3. This connection will play a central role in the present work and is a reason for the number theoretic character of Gaussian Gabor frames. The main references for the section are the Ramanujan notebooks edited by Berndt, in particular \cite[Chap.\ 17]{RamanujanIII} for the von Neumann lattice (theory of signature 2) and \cite[Chap.\ 33]{RamanujanV} for the hexagonal lattice (corresponding theory of signature 3). The part of the corresponding theory which we need has been put on solid ground by Borwein and Borwein \cite{BorBor_Cubic_91} (see also \cite{Borwein_AGM_1987}).

For what follows, we need to introduce Gauss' hypergeometric function ${}_{2}{F}_{1}$. For $n \in \Z$, let the rising Pochhammer symbol be denoted by
\begin{equation}
	(z)_n =\frac{\Gamma(z+n)}{\Gamma(z)},\ z \in \C,
\end{equation}
where $\Gamma(z)$ is Euler's Gamma function;
\begin{equation}
	\Gamma(z)=\int_{\mathbb{R}_{+}}t^{z-1}e^{-t}\ dt,\,\,\, \mbox{for}~ \mathrm{Re}(z)>0. 
\end{equation}
In this work, we will consider the case of Gauss' hypergeometric function with positive parameters $a,b,c$ and real variable $0<x<1$, defined by
\begin{equation}
	{}_{2}{F}_{1}(a,b; c; x)=\sum_{n=0}^{\infty} \frac{(a)_n(b)_n}{(c)_n}\frac{x^n}{n!}.
\end{equation}
The parameters fulfill $a+b \leq c$ and $x \in (0,1)$, which means that we do not run into convergence issues.
For the sake of completeness, and to justify the name \textit{elliptic modulus}, we introduce the complete elliptic integral of the first kind, $K$, and refer to the textbook of Whittaker and Watson \cite[Chap.\ 22.3]{WhiWat69}:
\begin{equation}\label{eq:K(k)}
	K(k) = \int_{0}^{\pi/2} \frac{d\varphi}{\sqrt{1-k^2\sin^2\varphi}}.
\end{equation}
Here, $0<k<1$ is the elliptic modulus of $K$. The elliptic modulus is also defined by means of the Jacobian elliptic functions, i.e., the theta constants, as
\begin{equation}
	k(q)=\frac{\theta_2(q)^2}{\theta_3(q)^2}, \quad 0 < |q| < 1.
\end{equation}
Obviously, it depends upon $q=e^{\pi i \tau},\ \tau \in \mathbb{H}$, but we will suppress this dependence in the sequel. The complementary elliptic modulus is denoted by $k'$ (and depends on $q$). It is defined by the property
\begin{equation}
	k^2+k'^2=1
	\quad \Longleftrightarrow \quad
	\frac{\theta_2(q)^4}{\theta_3(q)^4} + \frac{\theta_4(q)^4}{\theta_3(q)^4} = 1.
\end{equation}
The next formula is quite remarkable and was already known to Gauss
\begin{equation}\label{eq:K(k)_agm}
		\frac{2}{\pi} \, K(k) = \frac{2}{\pi} \, \int_0^{\pi/2} \frac{d \varphi}{\sqrt{1 - k^2 \sin(\varphi)^2}} = \sum_{n=0}^\infty \frac{\left( \tfrac{1}{2} \right)_n^2}{n!} \frac{k^{2n}}{n!} = {}_2F_1( \tfrac{1}{2}, \tfrac{1}{2}, 1, k^2) = \frac{1}{\agm(k',1)}.
\end{equation}
A proof, which is essentially the original proof of Gauss with some details filled in by Jacobi, is given in \cite[Chap.\ 1.2]{Borwein_AGM_1987} and \cite{Cox84}. The following result is not needed immediately for our purpose, but it will be essential in proving \eqref{eq:Landau_ag3}, i.e., that the conjectural value $\mathcal{L}_+$ of Landau's constant can be obtained as a cubic arithmetic-geometric mean (see \cite{Lan29} for the problem and \cite{Rad43} for the conjectural solution). The formula was known to Gauss and it can be found as Entry 34 in the textbook of Berndt \cite[Chap.\ 10]{RamanujanII};
\begin{equation}\label{eq:2F1_Gamma}
	{}_2F_1 \left(a,b; \tfrac{1}{2} (a+b+1); \tfrac{1}{2} \right) = \sqrt{\pi} \, \frac{\Gamma \left( \tfrac{a+b+1}{2} \right)}{\Gamma \left( \tfrac{a+1}{2} \right)\Gamma \left( \tfrac{b+1}{2} \right)}.
\end{equation}
Finally, we arrive at a kind of inversion formula, which was established by Ramanujan (see \cite[Entry 3, Chap.\ 17]{RamanujanIII}). Let $k = \theta_2^2/\theta_3^2$ be the elliptic modulus from above, then
\begin{equation}\label{eq:K_inversion}
	{}_2F_1 \left( \tfrac{1}{2},\tfrac{1}{2}; 1; k^2 \right) = \theta_3(q)^2.
\end{equation}
Note that this connects the theta constant to the arithmetic-geometric mean, and the complete elliptic integral of the first kind by \eqref{eq:K(k)_agm}.

The equation\eqref{eq:K_inversion} is a statement about lattice theta functions: $\theta_2(e^{-\pi \alpha})^2 = \theta_{\L_{1\times1}}(b_\circ; \alpha)$ and $b_\circ = (1/2, 1/2)$, yielding the deep hole in the von Neumann lattice, which is also the minimizer of the $\theta_{\L_{1\times1}}$ for any fixed $\alpha > 0$ (see \cite{BetKnu_Born_18}, \cite{FauSte17}, \cite{Jan96}). We remark that a deep hole is a point which globally maximizes the distance to the closest lattice point (local maximizers are called shallow holes) \cite[Chap.\ 1]{ConSlo98}. Evaluating $\theta_{\L_{1\times1}}$ at a lattice point, which always is a maximizer (see \cite{BetFauSte21}), we see that $\theta_3(e^{-\pi \alpha})^2 = \theta_{\L_{1\times1}}(0;\alpha)$. So, the elliptic modulus is the ratio
\begin{equation}
	k(e^{-\pi \alpha}) = \frac{\theta_{\L_{1\times1}}(b_\circ; \alpha)}{\theta_{\L_{1\times1}}(0;\alpha)},
\end{equation}
which is by \eqref{eq1:latticetheta} the ratio of the global minimum and maximum of the lattice theta function over a scaled von Neumann lattice. Gauss' hypergeometric function separates the minimum and maximum from each other.

Next, we show that the theory of signature 3 is intimately related to the lattice theta function over the hexagonal lattice $\L_2= H\Z^2$, where
\begin{equation*}
	H = \tfrac{\sqrt{2}}{\sqrt[4]{3}} \begin{pmatrix}
		1 & \frac{1}{2} \\
		0 & \frac{\sqrt{3}}{2}
	\end{pmatrix}
\end{equation*}
is the generating matrix for the lattice $\L_2$. It is pertinent to note that sometimes $\L_2$ is called a triangular lattice since half of its fundamental domain is an equilateral triangle, while the name  ``hexagonal lattice" is related to the fact that its Voronoi cells are all regular hexagons. The Voronoi cell of a lattice point is the set of all points which are closer to the lattice point than to any other point in the lattice (see \cite[Chap.\ 1]{ConSlo98}).

For $\alpha > 0$, $c = H \widetilde{c} = H(\widetilde{c}_1,\widetilde{c}_2)$, the hexagonal lattice theta function is explicitly given by
\begin{equation}\label{eq:hexa_latticesum}
	\theta_{\L_2}(c;\alpha) = \sum_{\l \in \L_2} e^{-\pi \alpha|\l + c|^2} = \sum_{(k,l) \in \Z^2} e^{\frac{-2\pi}{\sqrt{3}} ( (k+\widetilde{c}_1)^2 + (k+\widetilde{c}_1)(l+\widetilde{c}_2) + (l+\widetilde{c}_2)^2)}.
\end{equation}

Our next aim is to provide another significant relation between the hexagonal lattice theta function and the hypergeometric function, which will be obtained by Ramanujan's corresponding theory of signature 3. Ramanujan's corresponding theories of signature $r$ rely on
\begin{equation}
	{}_{2}{F}_{1}\big(\tfrac{1}{r},\tfrac{r-1}{r}; 1; . \big),
	\quad \text{ for } r \in \{2,3,4,6\}.
\end{equation}
This theory of signature 3 involves the cubic analogues of the squares of Jacobi's theta functions denoted as $a(q), b(q)$ and $c(q)$ as introduced in Section \ref{sec:analogues}. The cubic analogue of the elliptic modulus $k(q)$ and the complementary elliptic modulus $k'(q)$ are given by
\begin{equation}
	s(q) = \frac{c(q)}{a(q)}
	\quad \text{ and } \quad
	s'(q) = \frac{b(q)}{a(q)},
	\quad \text{ fulfilling } \quad
	s(q)^3 + s'(q)^3 = 1.
\end{equation}
The theory was put on a solid mathematical basement in \cite{BorBor_Cubic_91} and we also find the following formula \cite[Thm.\ 2.2(b) \& Thm.\ 2.3]{BorBor_Cubic_91} (see also \cite[Chap.\ 33]{RamanujanV}):
\begin{equation}\label{eq:signature3}
	{}_{2}{F}_{1}\Big(\tfrac{1}{3},\tfrac{2}{3}; 1; s^3\Big) = \frac{1}{\agc(1,s')} = a(q).
\end{equation}
The statement of \cite[Thm.\ 2.2(b)]{BorBor_Cubic_91} contains a harmless typo, which becomes apparent when looking at the proof. We note that \eqref{eq:signature3} is a statement on hexagonal lattice theta functions. Setting $q= e^{- \frac{2 \pi}{\sqrt{3}} \alpha}$, we have the following relations between the cubic analogues and hexagonal lattice theta functions:
\begin{equation}
	a(q) = \theta_{\L_2}(0;\alpha), \quad b(q) = \widehat{\theta}_{\L_2}(c_\circ;\alpha), \quad c(q) = \theta_{\L_2}(c_\circ,\alpha).
\end{equation}
The point $c_\circ = H(1/3,1/3)$ is a deep hole of the hexagonal lattice. Moreover, for all $\alpha >0$, the deep hole is the minimizer of the $\theta_{\L_2}$ and $\widehat{\theta}_{\L_2}$. This shows the intimate relation of the family of hexagonal lattice theta functions to Ramanujan's corresponding theory of signature~3.

\subsection{The constants \texorpdfstring{$G$}{G} and \texorpdfstring{$\mathcal{L}$}{L} and proof of Equation \texorpdfstring{\bf $\eqref{eq:Landau_ag3}$}{}}

In this section we prove \eqref{eq:Landau_ag3} from the \nameref{sec:remark}. We start with recalling that the conjectural value of Landau's constant $\mathcal{L}$, presented in \cite{Lan29}, is the following value found in \cite{Rad43};
\begin{equation}
	\mathcal{L}_+ = \frac{\Gamma(\frac{1}{3}) \Gamma(\frac{5}{6})}{\Gamma(\frac{1}{6})} \approx 0.543259 \ldots \, .
\end{equation}
It is known that
	$\frac{1}{2} < \mathcal{L} \leq \mathcal{L}_+$
and 
 conjectured that the second inequality is sharp.

\begin{proof}[Proof of \eqref{eq:Landau_ag3}]
	We use Gauss' formula \eqref{eq:2F1_Gamma}, connecting the hypergeometric function ${}_2F_1$ with the ratio of Gamma functions. The specific values we need are $a=1/3$ and $b=2/3$:
	\begin{equation}\label{eq:proof_2F1_aux}
		{}_2F_1(\tfrac{1}{3},\tfrac{2}{3};1;\tfrac{1}{2}) = \sqrt{\pi} \frac{\Gamma(1)}{\Gamma(\frac{2}{3}) \Gamma(\frac{5}{6})} = \frac{\Gamma(\frac{1}{2})}{\Gamma(\frac{2}{3}) \Gamma(\frac{5}{6})}.
	\end{equation}
	Now, we use Legendre's duplication formula for the Gamma function \cite[6.1.18]{AbrSte72}:
	\begin{equation}
		\Gamma(2z) = 2^{2z-1} \pi^{-1/2} \, \Gamma(z) \Gamma(z+\tfrac{1}{2}).
	\end{equation}
	Evaluating at $z= 1/6$ and some simple manipulations yield
	\begin{equation}
		\Gamma(\tfrac{2}{3}) = 2^{2/3} \frac{\Gamma(\frac{1}{2}) \Gamma(\frac{1}{3})}{\Gamma(\frac{1}{6})}.
	\end{equation}
	Substituting $\Gamma(2/3)$ in \eqref{eq:proof_2F1_aux} by the above expression then gives
	\begin{equation}
		{}_2F_1(\tfrac{1}{3},\tfrac{2}{3};1;\tfrac{1}{2}) = 2^{-2/3} \frac{\Gamma(\frac{1}{6})}{\Gamma(\frac{1}{3}) \Gamma(\frac{5}{6})} = \frac{1}{2^{2/3} \mathcal{L}_+}.
	\end{equation}
	By using the (symplectic) Poisson summation formula, we see that
		$b(e^{-\frac{2 \pi}{\sqrt{3}}}) = c(e^{-\frac{2 \pi}{\sqrt{3}}})$
	and by using the cubic identity, we find
	\begin{equation}
		a(e^{-\frac{2 \pi}{\sqrt{3}}})^3 = b(e^{-\frac{2 \pi}{\sqrt{3}}})^3 + c(e^{-\frac{2 \pi}{\sqrt{3}}})^3 = 2 b(e^{-\frac{2 \pi}{\sqrt{3}}})
		\quad \Longrightarrow \quad
		a(e^{-\frac{2 \pi}{\sqrt{3}}}) = 2^{1/3} b(e^{-\frac{2 \pi}{\sqrt{3}}}).
	\end{equation}
	Moreover, we have
	\begin{equation}
		s(e^{-\frac{2 \pi}{\sqrt{3}}})^3 = s'(e^{-\frac{2 \pi}{\sqrt{3}}})^3 = \frac{1}{2}.
	\end{equation}
	Plugging this into \eqref{eq:signature3} we obtain
	\begin{equation}
		{}_2F_1(\tfrac{1}{3},\tfrac{2}{3};1;\tfrac{1}{2}) = a(e^{-\frac{2 \pi}{\sqrt{3}}}) = 2^{-2/3} \frac{\Gamma(\frac{1}{6})}{\Gamma(\frac{1}{3}) \Gamma(\frac{5}{6})} = \frac{1}{2^{2/3} \mathcal{L}_+} = \frac{1}{\agc(1,2^{-1/3})}.
	\end{equation}
	So, we already arrived at the equality $2^{2/3} \mathcal{L}_+ = \agc(1, 2^{-1/3})$. Now, we use that $\agc$ is homogeneous and get
	\begin{equation}
		2 \mathcal{L}_+ = 2^{1/3} \agc(1, 2^{-1/3}) = \agc(2^{1/3}, 1).
	\end{equation}
	\vspace*{-0.82cm}
\end{proof}
We note that the connection $b(e^{-\frac{2\pi}{\sqrt{3}}}) = 1/(2 \mathcal{L}_+)$ was already established by one of the authors in \cite{Fau21-RJ} and we mainly filled in some details. Also, we want to briefly mention that Gauss' constant $G = \theta_4(e^{-\pi})^2$ appears in Landau's problem as well \cite{Fau21-RJ} (see \cite{Baernstein_Metric_2005}, \cite{Eremenko_Hyperbolic_2011} for the restricted problems in question). Originally, it appeared as the ratio of the arc length of the lemniscate of Bernoulli to the arc length of the unit circle. The computation of the arc length of the lemniscate involves elliptic integrals and Gauss computed it by hand:
\begin{equation}
	2 \varpi = 4 \int_0^{\pi/2} \frac{d x}{\sqrt{2 \left((1 - \frac{1}{2} \sin(x)^2\right)}} = 5.24412 \ldots \, .
\end{equation}
Thus, the ratio of $2 \varpi$ to the arc length of the unit circle is
	$\frac{2 \varpi}{2 \pi} = 0.834627 \ldots = G .$
Gauss observed that this is numerically $1/\agm(\sqrt{2},1)$. In fact, he was so enthusiastic about the connection he had found, that on May 30, 1799 he wrote in his diary (translated from \cite{Kle03}):

\medskip
\noindent
\textit{We have established that the arithmetic-geometric mean between 1 and $\sqrt{2}$ is $\frac{\pi}{\varpi}$ to the $11^{\text{th}}$ decimal place; the demonstration of this fact will surely open an entirely new field of analysis.}

\medskip
\noindent
We refer to the article by Cox \cite{Cox84} for more historical background and the proof Gauss actually found for his observation. Insights into the AGM for complex numbers are also given there.

\begin{figure}[h!]
	\hfill
	\subfigure[The lemniscate of Bernoulli]{
		\includegraphics[width=.4\textwidth]{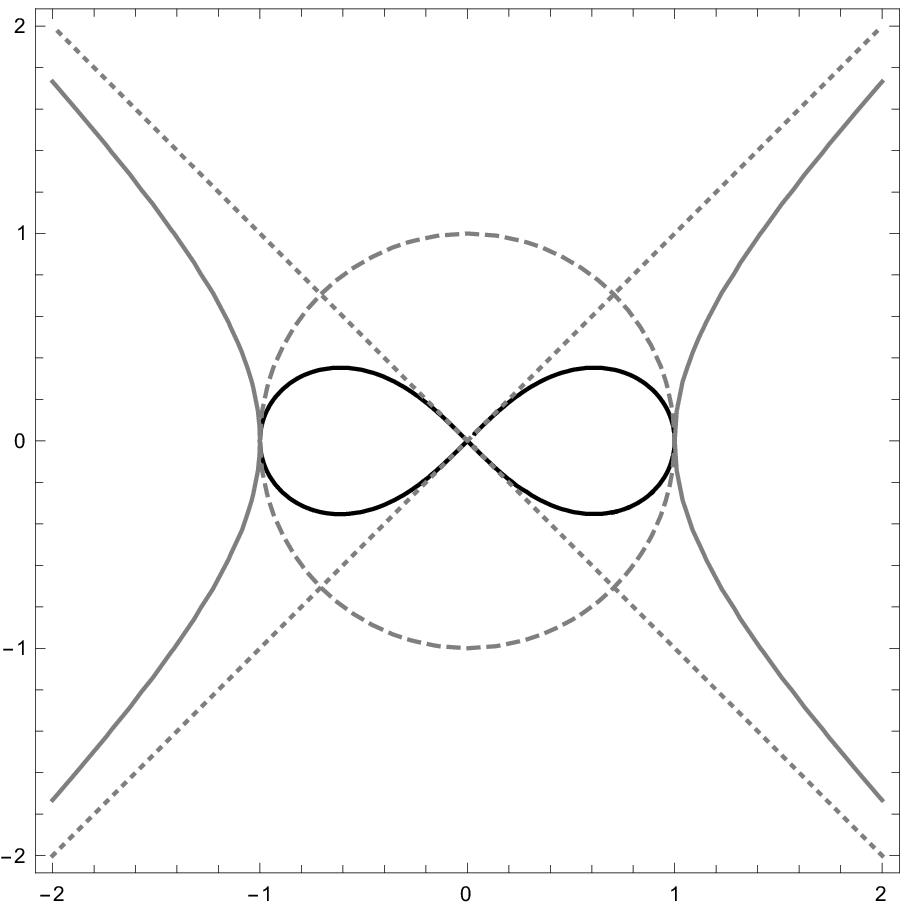}
	}
	\hfill
	\subfigure[Erdös lemniscate of degree 3]{
		\includegraphics[width=.4\textwidth]{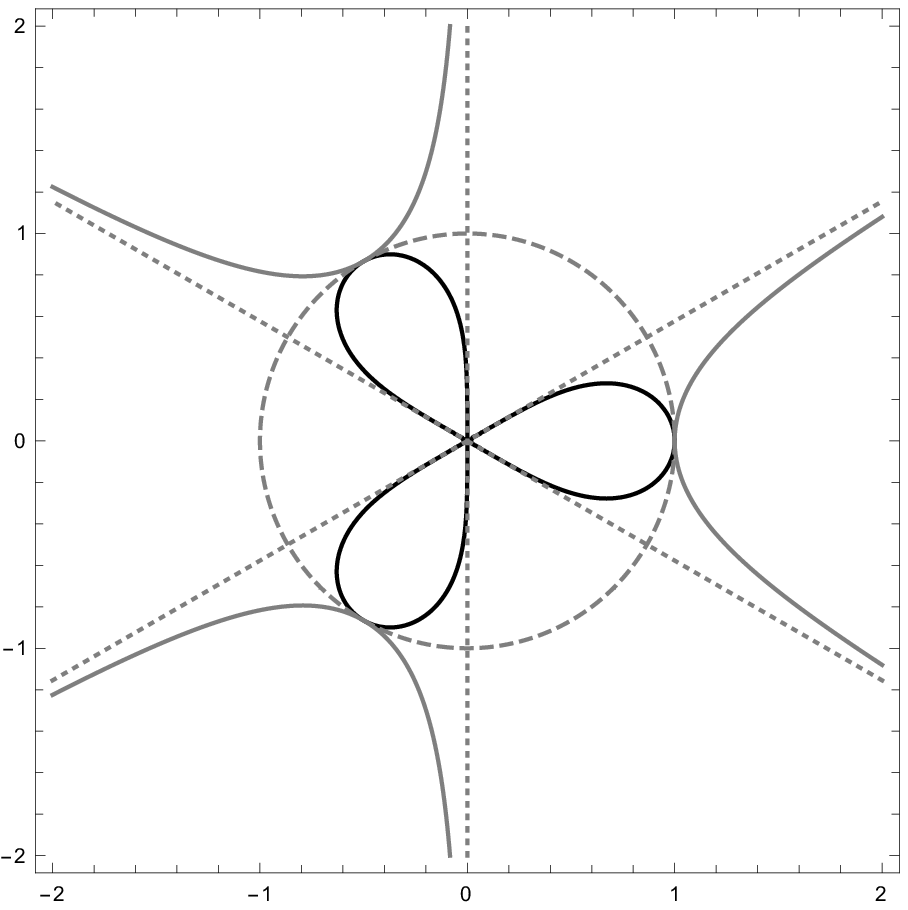}
	}
	\hspace*{\fill}
	\caption{The lemniscate of Bernoulli and Erdös lemniscate of degree 3 (put to scale). The reflection of the respective lemniscate with respect to the unit circle (inverting distances) yields the corresponding hyperbolas. The intersections of their asymptotes with the unit circle yield root systems. These are generating a (scaled and rotated) von Neumann lattice and a hexagonal lattice, respectively.}\label{fig:lemniscate}
\end{figure}

It stands to reason that the constant $\mathcal{L}_+$ will also play a role for the problem of finding the longest polynomial lemniscate of degree 3 (see Figure \ref{fig:lemniscate}). For a polynomial $p(z)$ of degree $n \geq 1$, the polynomial lemniscate is (see \cite{FryNaz09})
\begin{equation}
	L_p = \{z \in \C \mid |p(z)| = 1\}.
\end{equation}

An open conjecture of Erdös, Herzog, and Piranian \cite[Problem 12]{ErdHerPir58} states that for any fixed $n \geq 1$, the polynomial $p_0(z) = z^n - 1$ gives the maximal length of the lemniscate $L_p$ among all lemniscates of degree $n$. The lemniscates $L_{p_0}$ are also called Erdös lemniscates.

\section{Gabor Frame Bounds and the AGM}\label{sec:AGM}
\noindent
Our interest lies in computing Gaussian Gabor frame bounds. The results below show that they obey the AGM machinery and, hence, are intimately connected to the theory of elliptic integrals, hypergeometric functions and analytic number theory. However, this is only a side remark and we will not study these connections in depth here. We remark again that, classically, the frame operator is defined without the normalizing constant $\vol(\L)$. However, we want to define it in this way for two reasons:

\begin{enumerate}[(a)]
	\item In this way the frame bounds fulfill $A_\L \leq 1 \leq B_\L$, rather than $A_\L \leq \vol(\L)^{-1} \leq B_\L$.
	
	\item We do not need to re-normalize results, applicable to the frames bounds, which have been obtained by Gauss himself.
\end{enumerate}

\subsection{Computing Gabor frame bounds}
Determining whether the Gabor frame operator satisfies the spectral inequality \eqref{eq:frame} is a difficult task, at best. There have been several theoretical methods developed to investigate whether a Gabor system is a Gabor frame. In the case of integer or rational density, a popular method is the Zak transform \cite[Chap.\ 8]{Gro01}, which was introduced by Zak in the context of solid-state physics \cite{Zak67}. Alternatively, one can turn to duality theory. We present only the most essential results which are needed for our purpose. For more details we refer to the results of Daubechies, Landau and Landau \cite{DauLanLan94}, Janssen \cite{Jan95}, Ron and Shen \cite{RonShe97} and Wexler and Raz \cite{WexRaz90}. For a more recent treatise of the topic we refer to \cite{GroKop19}.

We have already seen how we can use the relaxed condition \eqref{eq:frame_relaxed} and FIGA \eqref{eq:FIGA} to estimate the spectral bounds $A_\L$ and $B_\L$ of the Gabor frame operator $S_\L$. We will now use the mentioned duality theory, which usually involves the (symplectic) Poisson summation formula, to compute exact bounds. By duality, we mean characterizing the Gabor systems $\G(g, \Lambda)$ which satisfy the frame inequality by using a condition on the system $\G(g,\Lambda^\circ)$.
%
Part of the theory allows us to relate the spectral bounds of the Gabor frame operator to sums over the adjoint lattice. We use a result of Janssen \cite{Jan96} in the generalized form given in \cite{Fau18} (see also \cite{FauSha22} where the method is described in some detail), tailored to our situation and normalization.
\begin{proposition}[Janssen \cite{Jan96}]\label{thm:Janssen_duality}
	For the Gaussian window $\varphi(t) = 2^{1/4} e^{-\pi t^2}$ and a lattice $\Lambda(\alpha) \subseteq \R^2$ of density $\alpha = 2N$, $N \in \N$, consider the Gabor system $\G(\varphi, \L)$. Then the spectral bounds of the associated frame operator $S_\L$ are given by
	\begin{align}\label{eq:Vgg}
		\begin{split}
			A_\L & = \min_{z \in \R^2} \sum_{\lambda^\circ\in\Lambda^\circ} e^{-\pi |\l^\circ|^2} e^{2 \pi i \sigma(\lambda^\circ,z)},\\
			B_\L & = \max_{z \in \R^2} \sum_{\lambda^\circ\in\Lambda^\circ} e^{-\pi |\l^\circ|^2} e^{2 \pi i \sigma(\lambda^\circ,z)}.
		\end{split}
	\end{align}
\end{proposition}
Above, the skew-symmetric form $\sigma(\l^\circ, z) = -\sigma(z, \l^\circ) = \l^\circ \cdot \mathcal{J} z$ appears as a consequence of the symplectic structure of the time-frequency plane (see \cite{Fol89}, \cite{Gos11}). However, this is merely a rotation of the co-ordinates and we could write the result using the dual lattice $\L^\perp$ and the Euclidean inner product $\l^\perp \cdot z$ in the complex exponential (see also \cite{Fau18}).

\subsection{The proof of Theorem \ref{thm:agm2}}
\begin{proof}
	For the proof, we use Proposition \ref{thm:Janssen_duality}. The problem is to find extremizers for
	\begin{equation}
		\widehat{\theta}_{\L_{1\times1}(1)}(b;n) = \sum_{k,l \in \Z} e^{-\pi n (k^2+l^2)} e^{2 \pi i (k b_2 - l b_1)}, \quad b=(b_1,b_2) \in \R^2.
	\end{equation}
	Note, that computing spectral bounds for lattice density $2N$, $N \in \N$, corresponds to choosing positive integer values for the parameter $\alpha$ of the lattice theta function. The discrepancy of the factor 2 mainly comes from the fact that
	\begin{equation}
		|\langle \pi(z) \varphi, \pi(\l) \varphi \rangle| = e^{-\frac{\pi}{2} |\l-z|^2}.
	\end{equation}
	
	An application of the triangle inequality shows  that the maximum of $\widehat{\theta}_\L (b;\alpha)$ (for any $\L$ and any $\alpha > 0$) is always attained for $b = (0,0)$ and, by periodicity at any other lattice point. Determining the exact value of $b$ which gives the minimum is in general a challenging task \cite{BetFauSte21} (see also \cite{Baernstein_HeatKernel_1997}). Due to the occurring symmetries in the von Neumann lattice, the minimum is always attained in the center of a fundamental square, i.e., at $b \in (1/2,1/2)+\Z^2$ (compare \cite{Jan96} and also \cite{BetKnu_Born_18}) times the right scaling factor. For the scaled von Neumann lattice $\L_{1\times1}(\alpha) = {\alpha}^{-1/2} \Z \times {\alpha}^{-1/2} \Z$, $\alpha = 2N$ with $N \in \N$, we obtain the explicit Gaussian Gabor frame bounds, where we set $q = e^{-\pi}$:
	\begin{align}
		A_{1\times1}(\alpha) & = \sum_{k,l \in \Z} (-1)^{k+l} e^{-\pi N (k^2+l^2)} = \theta_4(q^N)^2,\\
		B_{1\times1}(\alpha) & = \sum_{k,l \in \Z} e^{-\pi N (k^2+l^2)} = \theta_3(q^N)^2.
	\end{align}
	Now, if $\alpha$ is a (positive) power of 2, say $\alpha = 2^n$, $n \in \N_+$, then we have Gabor systems over the following scaled family of von Neumann lattices:
	\begin{equation}
		\L_{1\times1}(2^{n}) = 2^{-n/2} \Z \times 2^{-n/2} \Z.
	\end{equation}
	We obtain the spectral bounds of the Gabor frame operator as
	\begin{align}
		A_{1\times1}(2^n) = \theta_4(q^{2^n})^2
		\quad \text{ and } \quad
		B_{1\times1}(2^n) = \theta_3(q^{2^n})^2.
	\end{align}
	All that is left to observe is that we can now use the arithmetic-geometric mean iteration for Jacobi's theta functions:
	\begin{equation}
		\theta_3(q)^2+\theta_4(q)^2 = 2 \theta_3(q^2)^2
		\quad \text{ and } \quad
		\theta_3(q) \theta_4(q) = \theta_4(q^2)^2.
	\end{equation}
	This implies that
	\begin{equation}
		B_{1\times1}(2^{n+1}) = \frac{A_{1\times1}(2^n)+B_{1\times1}(2^n)}{2}
		\quad \text{ and } \quad
		A_{1\times1}(2^{n+1}) = \sqrt{A_{1\times1}(2^n) B_{1\times1}(2^n)}.
	\end{equation}
\end{proof}

The spectral bounds for the frame operator for $\L_{1\times1}(2) = 2^{-1/2} \Z \times 2 ^{-1/2} \Z$ have a simple algebraic dependence: we know that (see \cite{Faulhuber_Hexagonal_2018, StrBea03})
\begin{equation}
	\kappa_{1\times1}(2) = \frac{B_{1\times1}(2)}{A_{1\times1}(2)} = \sqrt{2}
	\quad \Longleftrightarrow \quad
	B_{1\times1}(2) = \sqrt{2} A_{1\times1}(2).
\end{equation}
This goes back to the Jacobi identity $\theta_2(q)^4+\theta_4(q)^4 = \theta_3(q)^4$.
We can now use our knowledge that $A_{1\times1}(2^n)$ and $B_{1\times1}(2^n)$ obey the $\agm$ machinery to iteratively compute condition numbers $\kappa_{1\times1}(2^n) = B_{1\times1}(2^n)/A_{1\times1}(2^n)$ of the associated family of frame operators.
\begin{equation}
	B_{1\times1}(2^{n+1}) = (1+\kappa_{1\times1}(2^n)) \frac{A_{1\times1}(2^n)}{2}
	\quad \text{ and } \quad
	A_{1\times1}(2^{n+1}) = \sqrt{\kappa_{1\times1}(2^n)} \, A_{1\times1}(2^n).
\end{equation}
It readily follows that the sequence of condition numbers obeys the iterative rule
\begin{equation}
	\kappa_{1\times1}(2^{n+1}) = \frac{1}{2} \left( \frac{1}{\sqrt{\kappa_{1\times1}(2^n)}} + \sqrt{\kappa_{1\times1}(2^n)} \right).
\end{equation}
As $\kappa_{1\times1}(2^n) \geq 1$, it is simple to see that the sequence of condition numbers is decreasing:
\begin{equation}
	\kappa_{1\times1}(2^n) \geq \sqrt{\kappa_{1\times1}(2^n)} \geq \frac{1}{2} \left( \frac{1}{\sqrt{\kappa_{1\times1}(2^n)}} + \sqrt{\kappa_{1\times1}(2^n)} \right) = \kappa_{1\times1}(2^{n+1}).
\end{equation}
In fact, it is converging to 1. This shows that the frame operator converges to the identity operator as the density tends to infinity. We note that this is only a special case of a result of Feichtinger and Zimmermann \cite{FeiZim98}. However, the above proof has already been accessible to Gauss and implies the existence of Gabor frames if the density of the lattice is high enough. Obviously, it considerably precedes the theory of Gabor frames, which makes it especially interesting. We will now recall a result of Gauss\footnote{From \cite{Gauss_WerkeIII}: \textit{Das Arithmetisch Geometrische Mittel zwischen $(Px)^2$ und $(Qx)^2$ ist allemal $= 1$}.} \cite[Theorem (23), p.467]{Gauss_WerkeIII} (see also \cite[Entry 102]{Kle03}) formulated in the language of Jacobi theta functions.
\begin{theorem*}[Gauss, 1799]
	The arithmetic-geometric mean of $\theta_3(q)^2$ and $\theta_4(q)^2$ always equals~1.
\end{theorem*}

Note that Gauss' theorem shows that the spectral bounds $A_{1\times1}(2^n)$ and $B_{1\times1}(2^n)$ of the frame operator $S_{\L_{1\times1}}$ both tend to 1 as $n \to \infty$. Hence, the double-sided operator inequality
\begin{equation}
	A_n \mathsf{I}_{L^2} \leq S_{\L_{2^{-n}}} \leq B_n \mathsf{I}_{L^2}
\end{equation}
shows that $S_{\L_{2^{-n}}} \to \mathsf{I}_{L^2}$ in the operator norm. We know that the rate of convergence is quadratic \cite{BorBor_Cubic_91}, which fits nicely with doubling the lattice density in each iteration of $\agm$:
\begin{equation}
	\agm(B_{1\times1}(2), A_{1\times1}(2)) = \agm(B_{1\times1}(2^n), A_{1\times1}(2^n)) = 1.
\end{equation}

\subsection{The proof of Theorem \ref{thm:agm3}}

\begin{proof}
	We begin with a short overview of the hexagonal lattice. It can be realized as
	\begin{equation}
		\LHex (\alpha) = \alpha^{-1/2} H\Z^2,
		\quad \text{ where } \quad
		H = \tfrac{\sqrt{2}}{\sqrt[4]{3}}
		\begin{pmatrix}
			1 & \frac{1}{2} \\
			0 & \frac{\sqrt{3}}{2}
		\end{pmatrix}, \, \alpha > 0.
	\end{equation}
	The circumcenter of the fundamental triangle with vertices $(0,0)$, $H(1,0)$ and $H(0,1)^T$ is easily verifiable to be the point $c_\circ = H(1/3,1/3) = 2^{-1/2}3^{-1/4}(1,\ 3^{-1/2})$.  The optimal frame bounds are determined with Janssen's result again. As mentioned ealier, for the scaled hexagonal lattice of density $2N$ for some $N\in\N$, it suffices to evaluate at the origin for the upper bound. Evaluating at the circumcenter of a fundamental triangle, for instance, at $\frac{1}{2N} c_\circ$, gives the lower frame due to a result of Baernstein \cite{Baernstein_HeatKernel_1997}. We obtain explicit Gaussian Gabor frame bounds. Setting $\alpha = 2N$ and $q = e^{-\frac{2 \pi}{\sqrt{3}}}$, we get (cf. Section \ref{sec:Ramanujan})
	\begin{align}
		A_2(\alpha) & = \sum_{k,l \in \Z} e^{2\pi i/3(k-l)} e^{-2\pi N/\sqrt{3}(k^2+kl+l^2)} = b(q^N),\\
		B_2(\alpha) & = \sum_{k,l \in \Z} e^{-2\pi N/\sqrt{3}(k^2+kl+l^2)} = a(q^N).
	\end{align}
	Now, if $\alpha$ is twice a (non-negative) power of 3, say $\alpha = 2\cdot 3^{n-1}$, $n \in \N$, then (similarly to the family of von Neumann lattices) we concern ourselves with the following family of lattices:  $\LHex(2\cdot 3^{n-1}) = 2^{-1/2}3^{-(n-1)/2} H\Z^2$. We see that the spectral bounds can be written as
	\begin{align}
		A_2 (2 \cdot 3^{n-1}) = b(q^{3^{n-1}})
		\quad \text{ and } \quad
		B_2 (2 \cdot 3^{n-1}) = a(q^{3^{n-1}}).
	\end{align}
	As established in Section \ref{sec:AGM_theta}, \eqref{eq:correspondence}, the sequences $\left(A_2(2\cdot 3^{n-1})\right)_{n\in\N}$ and $\left(B_2(2\cdot 3^{n-1})\right)_{n\in\N}$ comply to the cubic arithmetic-geometric mean construction;
	\begin{align}
		B_2(2 \cdot 3^n) & = \frac{B_2(2 \cdot 3^{n-1}) + 2 A_2(2 \cdot 3^{n-1})}{3}\\
		A_2(2 \cdot 3^n) & = \sqrt[3]{A_2(2 \cdot 3^{n-1}) \, \frac{A_2(2 \cdot 3^{n-1})^2 + A_2(2 \cdot 3^{n-1}) B_2(2 \cdot 3^{n-1}) + B_2(2 \cdot 3^{n-1})^2}{3}} \, .
	\end{align}
\end{proof}
A simple calculation with the Poisson summation formula relates $b(q)$ and $c(q)$ analogously to the famous formula for $\theta_4(q)$ and $\theta_2(q)$. In particular, it shows $B_2 (2)^3 = 2 A_2(2)^3 $. Equivalently, the condition number is $\kappa_2(2) = \sqrt[3]{2}$.
Inductively, we obtain
\begin{align}
	\kappa_2(2\cdot 3^{n})^3 =& \frac{1}{9}\cdot \frac{(B_2(2\cdot 3^{n-1})+2A_2(2\cdot 3^{n-1}))^3}{A_2(2\cdot 3^{n-1})(A_2(2\cdot 3^{n-1})^2+A_2(2\cdot 3^{n-1})B_2(2\cdot 3^{n-1})+B_2(2\cdot 3^{n-1})^2)} \\
	&= \frac{1}{9}\cdot \left(\tfrac{B_2(2\cdot 3^{n-1})}{A_2(2\cdot 3^{n-1})}+2\right)^3 \left(1+\tfrac{B_2(2\cdot 3^{n-1})}{A_2(2\cdot 3^{n-1})}+\left(\tfrac{B_2(2\cdot 3^{n-1})}{A_2(2\cdot 3^{n-1})}\right)^2\right)^{-1}\\
	=&\frac{1}{9}\cdot\left(\kappa_2(2\cdot 3^{n-1})+2\right)^3 \left(1+\kappa_2(2\cdot 3^{n-1})+\kappa_2(2\cdot 3^{n-1})^2\right)^{-1}, \\
	\kappa_2(2\cdot 3^{n}) =& \frac{\sqrt[3]{3}}{3}\cdot \left(\kappa_2(2\cdot 3^{n-1})+2\right) \left(1+\kappa_2(2\cdot 3^{n-1})+\kappa_2(2\cdot 3^{n-1})^2\right)^{-1/3}.
\end{align}

We observe that the sequence of condition numbers is decreasing because
\begin{equation}
	\Big(\frac{\kappa_2(2\cdot 3^{n})}{\kappa_2(2\cdot 3^{n-1})}\Big)^3 = \frac{1}{9}\cdot\frac{1+2 \kappa_2(2\cdot 3^{n-1})^{-1}}{\left(1+\kappa_2(2\cdot 3^{n-1})+\kappa_2(2\cdot 3^{n-1})^2\right)^{1/3}} < 
	\frac{1}{9}\cdot \frac{1+\tfrac{2}{1}}{\left(1+1+1\right)^{1/3}}
	<1.
\end{equation}
Independently of the result of Feichtinger and Zimmermann \cite{FeiZim98}, we will conclude that the frame operator tends to the identity. This follows from properties of cubic theta functions. 
\begin{theorem*}[Borwein and Borwein, 1991]
	The cubic arithmetic-geometric mean of $a(q)$ and $b(q)$ always equals 1.
\end{theorem*}
The proof can be found in \cite[Thm.\ 2.3]{BorBor_Cubic_91}. As for the von Neumann lattice, this implies the convergence $S_{\LHex(2\cdot 3^{n})} \to \mathsf{I}_{L^2}$, $n \to \infty$.
The rate of convergence is cubic \cite[Chap.~1.1]{Borwein_AGM_1987}, \cite{BorBor_Cubic_91}, which fits nicely with the triplication of the lattice density in each iteration of $\agc$:
\begin{equation}
	\agc(B_2(2), A_2(2)) = \agc(B_2(2 \cdot 3^n), A_2(2 \cdot 3^n)) = 1.
\end{equation}

\section{Rectangular Lattices}

\noindent
We will now present some results for rectangular lattices, as their frame bounds can be expressed by means of the frame bounds for the von Neumann lattice. We will handle lattices of even density.  In that case, the optimal frame bounds for a lattice $\L_{a \times b}(\alpha) = \alpha^{-1/2} (a \Z \times b \Z)$ with $\alpha = 2N$, $N \in \N$, $ab = 1$ and $q = e^{-\pi}$ are given by (see \cite{FauSte17})
\begin{align}
		A_{a \times b} (2N) & = \theta_4(q^{N a^2})\theta_4(q^{N b^2}),\\
		B_{a \times b} (2N)& = \theta_3(q^{Na^2})\theta_3(q^{Nb^2}).
	\end{align}
We observe what happens by doubling the density while preserving the initial ratio $a:b$.
	\begin{align}
		A_{a \times b} (2^{n+1})
		& = \theta_4(q^{2^{n} a^2})\theta_4(q^{2^{n} b^2})\\
		&= \sqrt{\theta_4(q^{2^{n-1} a^2}) \, \theta_3(q^{2^{n-1} a^2}) \, \theta_4(q^{2^{n-1} b^2}) \, \theta_3(q^{2^{n-1} b^2})}\\
		& = \sqrt{A_{a \times b} (2^n) \, B_{a \times b} (2^n)} \, .
\end{align}
The situation with the upper bound is a bit more involved.
	\begin{align}
		B_{a \times b} (2^{n+1}) & = \theta_3(q^{2^{n} a^2}) \theta_3(q^{2^{n} b^2})\\
		& = \left(\frac{\theta_3(q^{2^{n-1} a^2})^2+\theta_4(q^{2^{n-1} a^2})^2}{2}
		\, \frac{\theta_3(q^{2^{n-1} b^2})^2+\theta_4(q^{2^{n-1} b^2})^2}{2}\right)^{1/2} \\
		&=\tfrac{1}{2}\bigg(\theta_3(q^{2^{n-1} a^2})^2\theta_3(q^{2^{n-1} b^2})^2 + \theta_3(q^{2^{n-1} a^2})^2\theta_4(q^{2^{n-1} b^2})^2 \\
		&\qquad\qquad+ \theta_4(q^{2^{n-1} a^2})^2\theta_3(q^{2^{n-1} b^2})^2 
+ \theta_4(q^{2^{n-1} a^2})^2\theta_4(q^{2^{n-1} b^2})^2
		\bigg)^{1/2}.
\end{align}
As all standard relations between theta functions rely on an integer ratio between $a$ and $b$, further meaningful transformations are not readily available. We may simplify the quantities involving theta functions of the same kind:
\begin{equation}
	\theta_3(q^{2^{n-1} a^2})^2\theta_3(q^{2^{n-1} b^2})^2 = B_{a \times b}(2^n)^2
	\quad \text{ and } \quad
	\theta_4(q^{2^{n-1} a^2})^2\theta_4(q^{2^{n-1} b^2})^2 = A_{a \times b}(2^n)^2.
\end{equation}
However, it is unclear how to treat the terms involving theta functions of different kinds mixed together, as the lattice constants $a$ and $b$ only need to fulfill $a,b \in \R_+$, $ab = 1$.

\section{Two conjectures}\label{sec:conjecture}
\noindent
To our knowledge, the following conjecture has not appeared in a printed form before, but, nonetheless, the authors cannot claim any credit for coming up with this conjecture. It was brought forward to one of the authors in a private correspondence with Y.~Lyubarskii in 2018 when both were affiliated with NTNU Trondheim. Y.~Lyubarskii, on the other hand, was informed about the conjecture by T.~Strohmer already in 2009 in a private correspondence. The following conjecture should therefore, to the best of our knowledge, be addressed to T.~Strohmer from UC Davis and we publish it with the consent of the originator.
\begin{conjecture}[Strohmer, 2009]\label{con:Strohmer}
	Let $\kappa_3(\alpha)$ and $\kappa_4(\alpha)$ denote the frame condition of the Gaussian Gabor system with hexagonal $\L_2(\alpha) = \alpha^{-1/2} H \Z^2$ and von Neumann lattice $\L_{1\times1}(\alpha) = \alpha^{-1/2} \Z^2$, respectively, of density $\alpha$. Then, the result in \cite{BorGroLyu10} implies that
	\begin{equation}
		\kappa_3(\alpha) - \frac{C_3}{1-\frac{1}{\alpha}} \to 0
		\quad \text{ and } \quad
		\kappa_4(\alpha) - \frac{C_4}{1-\frac{1}{\alpha}} \to 0,
	\end{equation}
	as we approach the critical density $\alpha \to 1$ from above. The conjectural exact constants are
	\begin{equation}
		C_r = \frac{2 \pi r}{\tan(\frac{\pi}{r})} \frac{\Gamma(\frac{2}{r})^2}{\Gamma(\frac{1}{r})^4}, \quad r = 3,4.
	\end{equation}
\end{conjecture}
Plugging in the values, we obtain
\begin{equation}
	C_3 = \frac{6 \pi}{\sqrt{3}} \, \frac{\Gamma(2/3)^2}{\Gamma(1/3)^4} \approx 0.387438 \ldots
	\quad \text{ and } \quad
	C_4 = 8 \pi \, \frac{\Gamma(1/2)^2}{\Gamma(1/4)^4} \approx 0.456947 \ldots \, .
\end{equation}
This would provide evidence that the conjecture of Strohmer and Beaver raised in \cite{StrBea03} holds close to the critical density. At least it shows that asymptotically Gaussian Gabor frames over the hexagonal lattice have smaller condition numbers than over the square lattice:

Assume $\kappa_3(\alpha) < \kappa_4(\alpha)$ for all $0< \alpha < 1$ (necessary for the conjecture Strohmer and Beaver in \cite{StrBea03}). Then the asymptotics in Conjecture \ref{con:Strohmer}, proved in \cite{BorGroLyu10}, imply that for every $\varepsilon > 0$ there exists $\alpha_0 > 1$ such that for all $\alpha \in(1, \alpha_0)$ we have
\begin{equation}
	\abs{\kappa_3(\alpha) - \frac{C_3}{1-\frac{1}{\alpha}}} < \varepsilon
	\quad \text{ and } \quad
	\abs{ \kappa_4(\alpha) - \frac{C_4}{1-\frac{1}{\alpha}}} < \varepsilon.
\end{equation}
Rearranging the inequalities for $\kappa_3(\alpha)$ and $\kappa_4(\alpha)$ yields
\begin{equation}
	- \varepsilon + \frac{C_3}{1-\frac{1}{\alpha}} < \kappa_3(\alpha) < \varepsilon + \frac{C_3}{1-\frac{1}{\alpha}}
	\quad \text{ and } \quad
	- \varepsilon + \frac{C_4}{1-\frac{1}{\alpha}} <  \kappa_4(\alpha) < \varepsilon + \frac{C_4}{1-\frac{1}{\alpha}}.
\end{equation}
We obtain the following chain of inequalities:
\begin{equation}
	0 < \kappa_4(\alpha) - \kappa_3(\alpha) < \varepsilon + \frac{C_4}{1-\frac{1}{\alpha}} + \varepsilon - \frac{C_3}{1-\frac{1}{\alpha}} = 2 \varepsilon + \frac{C_4 - C_3}{1-\frac{1}{\alpha}}.
\end{equation}
As $\varepsilon > 0$ and $1-\frac{1}{\alpha} > 0$, this leads to
\begin{equation}
	C_4 + 2 \varepsilon (1-\tfrac{1}{\alpha}) > C_3
	\quad \Longrightarrow \quad
	C_4 > C_3, \; \text{ as } \varepsilon \to 0 \text{ and for } \alpha \to 1 \text{ from above}.
\end{equation}
So, we have shown that $C_4$ being larger than $C_3$ is a necessary condition for the conjecture of Strohmer and Beaver (as phrased in \cite{AbrDoe12}) to hold. In \cite{BorGroLyu10} it was also shown that the asymptotic behavior holds for rectangular lattices and there is no evidence that this behavior does not hold for arbitrary lattices. However, we are not aware of a proof for the general statement and consider it as an interesting problem. Still, we boldly formulate the following conjecture, based on the results of Montgomery \cite{Mon88} and Bétermin, Faulhuber, and Steinerberger \cite{BetFauSte21}.
\begin{conjecture}\label{con:kappa}
	Let $\kappa_\L(\alpha)$ be the condition number of the frame operator associated to the Gaussian Gabor system $\G(\varphi, \L(\alpha))$, where $\L(\alpha) \subset \R^2$ is a lattice of density $\alpha > 1$. Then, there exists a constant $C_\L$, depending on the geometry of the lattice $\L(1)$, such that
	\begin{equation}
		\kappa_\L(\alpha) - \frac{C_\L}{1-\frac{1}{\alpha}} \to 0,
	\end{equation}
	as we approach the critical density $\alpha \to 1$ from above. Furthermore, we conjecture that
	\begin{equation}\label{eq:conjecture}
		C_3 \leq C_\L,
	\end{equation}
	with equality if and only if $\L$ is the hexagonal lattice.
\end{conjecture}
\noindent
At the moment, it is not clear to the authors how one should approach Conjecture \ref{con:Strohmer} or \eqref{eq:conjecture} in Conjecture \ref{con:kappa}. One could use a Zak transform approach (or an equivalent approach) computing sequences of condition numbers $\kappa_3(1+1/n)$ and $\kappa_4(1+1/n)$, $n \in \N$. This will become very computationally intensive: we will need to compute and optimize eigenvalues of $n \times n$ matrices where the entries will be lattice theta functions. Maybe there is an approach via modular functions and the theories of Gauss, Ramanujan, and Borwein and Borwein.

%
%
%
%
\bibliographystyle{plain}

\begin{thebibliography}{10}

\bibitem{AbrSte72}
Milton Abramowitz and Irene Stegun.
\newblock{\em Handbook of Mathematical Functions}.
\newblock{Applied Mathematics Series 55}, 10$^{\text{th}}$ edition. National Bureau of Standards, U.S.\ Department of Commerce, 1972.

\bibitem{AbrDoe12}
L.\ Daniel Abreu and Monika Dörfler.
\newblock {An inverse problem for localization operators}.
\newblock {\em Inverse Problems}, 28(11):115001, 2012.
\href{https://doi.org/10.1088/0266-5611/28/11/115001}{10.1088/0266-5611/28/11/115001}

\bibitem{Baernstein_HeatKernel_1997}
Albert {Baernstein II}.
\newblock {A minimum problem for heat kernels of flat tori}.
\newblock In {\em {Extremal {R}iemann surfaces ({S}an {F}rancisco, {CA}, 1995)}}, volume 201 of {\em {Contemporary Mathematics}}, page 227--243. American Mathematical Society, Providence, RI, 1997.
\href{https://doi.org/10.1090/conm/201}{10.1090/conm/201}

\bibitem{Baernstein_Metric_2005}
Albert {Baernstein II}, Alexandre Eremenko, Alexander Fryntov, and Alexander Solynin.
\newblock {Sharp estimates for hyperbolic metrics and covering theorems of Landau type}.
\newblock {\em Annales Academiae Scientiarum Fennicae Mathematica}, 30:113–133, 2005.
\url{https://www.acadsci.fi/mathematica/Vol30/baernst.html}

\bibitem{Bal81}
Roger Balian.
\newblock {Un principe d’incertitude fort en théorie du signal ou en
  mécanique quantique}.
\newblock {\em CR Acad. Sci. Paris}, 292(2):1357–1361, 1981.
\url{http://gallica.bnf.fr/ark:/12148/bpt6k5657308b/f379.item}

\bibitem{BarButGirKla71}
Valentine Bargmann, Paolo Butera, Luciano Girardello, and John Klauder.
\newblock {On the completeness of the coherent states}.
\newblock {\em Reports on Mathematical Physics}, 2(4):221–228, 1971.
\href{https://doi.org/10.1016/0034-4877(71)90006-1}{10.1016/0034-4877(71)90006-1}

\bibitem{RamanujanII}
Bruce~C.\ Berndt.
\newblock {\em {Ramanujan’s Notebooks, Part II}}.
\newblock Springer, 1989.
\href{https://doi.org/10.1007/978-1-4612-4530-8}{10.1007/978-1-4612-4530-8}

\bibitem{RamanujanIII}
Bruce~C.\ Berndt.
\newblock {\em {Ramanujan’s Notebooks, Part III}}.
\newblock Springer, 1991.
\href{https://doi.org/10.1007/978-1-4612-0965-2}{10.1007/978-1-4612-0965-2}

\bibitem{RamanujanV}
Bruce~C.\ Berndt.
\newblock {\em {Ramanujan’s Notebooks, Part V}}.
\newblock Springer, 1998.
\href{https://doi.org/10.1007/978-1-4612-1624-7}{10.1007/978-1-4612-1624-7}

\bibitem{BetFau21}
Laurent Bétermin and Markus Faulhuber.
\newblock {Maximal Theta Functions – Universal Optimality of the Hexagonal Lattice for Madelung-Like Lattice Energies}.
\newblock {\em {Journal d'Analyse Math\'ematique}}, accepted, 2021.
\href{https://arxiv.org/abs/2007.15977}{arXiv:2007.15977}.

\bibitem{BetFauSte21}
Laurent Bétermin, Markus Faulhuber, and Stefan Steinerberger.
\newblock {A variational principle for Gaussian lattice sums}.
\newblock {\em arXiv preprint}, 2021.
\href{https://arxiv.org/abs/2110.06008}{arXiv:2110.060080}

\bibitem{BetKnu_Born_18}
Laurent Bétermin and Hans Knüpfer.
\newblock {On Born's conjecture about optimal distribution of charges for an infinite ionic crystal}.
\newblock {\em Journal of Nonlinear Science}, 28(5):1629–1656, 2018.
\href{https://doi.org/10.1007/s00332-018-9460-3}{10.1007/s00332-018-9460-3}

\bibitem{BooZakZuc83}
Michael Boon, Joshua Zak, and I.~John Zucker.
\newblock {Rational von Neumann lattices}.
\newblock {\em Journal of Mathematical Physics}, 24(2):316–323, 1983.
\href{https://doi.org/10.1063/1.525682}{10.1063/1.525682}

\bibitem{BorGroLyu10}
Alexander Borichev, Karlheinz Gröchenig, and Yurii Lyubarskii.
\newblock {Frame constants of Gabor frames near the critical density}.
\newblock {\em Journal de Mathématiques Pures et Appliqués}, 94(2):170--182, 2010.\\
\href{https://doi.org/10.1016/j.matpur.2010.01.001}{10.1016/j.matpur.2010.01.0010}

\bibitem{Borwein_AGM_1987}
Jonathan~M.\ Borwein and Peter~B.\ Borwein.
\newblock {\em {Pi and the AGM}}.
\newblock Wiley, New York, 1987.

\bibitem{BorBor_Cubic_91}
Jonathan~M.\ Borwein and Peter~B.\ Borwein.
\newblock {A Cubic Counterpart of Jacobi's Identity and the AGM}.
\newblock {\em Transactions of the American Mathematical Society}, 332(2):691–701, February 1991.
\href{https://doi.org/10.2307/2001551}{10.2307/2001551}

\bibitem{BorBorGar95}
Jonathan~M.\ Borwein, Peter~B.\ Borwein, and Frank G.\ Garvan.
\newblock{Some Cubic Modular Identities of Ramanujan}
\newblock{\em Transactions of the American Mathematical Society}, 343(1):35--47, 1995.
\href{https://doi.org/10.2307/2154520}{10.2307/2154520}

\bibitem{CohKum07}
Henry Cohn and Abhinav Kumar.
\newblock {Universally optimal distribution of points on spheres}.
\newblock {\em Journal of the American Mathematical Society}, 20(1):99--148, 2007.
\href{https://doi.org/10.1090/S0894-0347-06-00546-7}{10.1090/S0894-0347-06-00546-7}

\bibitem{ConSlo98}
John~H. Conway and Neil J.~A.\ Sloane.
\newblock {\em {Sphere Packings, Lattices and Groups}}, volume 290 of
{\em {Grundlehren der Mathematischen Wissenschaften [Fundamental Principles of Mathematical Sciences]}}.
\newblock Springer-Verlag, New York, 3. edition, 1998.
\href{https://doi.org/10.1007/978-1-4757-6568-7}{10.1007/978-1-4757-6568-7}

\bibitem{Cox84}
David~A.\ Cox.
\newblock {The arithmetic-geometric mean of Gauss}.
\newblock {\em {L'Enseignement Mathématique}}, 30:275–330, 1984.
\href{https://doi.org/10.1007/978-1-4757-3240-5_55}{10.1007/978-1-4757-3240-5\_55}

\bibitem{DauLanLan94}
Ingrid Daubechies, Henry~J.\ Landau, and Zeph Landau.
\newblock {Gabor Time-Frequency Lattices and the Wexler–Raz Identity}.
\newblock {\em Journal of Fourier Analysis and Applications}, 1(4):437--478, 1994.\\
\href{https://doi.org/10.1007/s00041-001-4018-3}{10.1007/s00041-001-4018-30}

\bibitem{ErdHerPir58}
Paul Erdös, Fritz Herzog, and George Piranian.
\newblock {Metric properties of polynomials}.
\newblock {\em Journal d'Analyse Mathematique}, 6(1):125–148, 1958.
\href{https://doi.org/10.1007/BF02790232}{10.1007/BF02790232}

\bibitem{Eremenko_Hyperbolic_2011}
Alexandre Eremenko.
\newblock {On the hyperbolic metric of the complement of a rectangular lattice}.
\newblock {\em arXiv preprint}, 2011.
\href{https://arxiv.org/abs/1110.2696}{arXiv:1110.2696}

\bibitem{Fau18}
Markus Faulhuber.
\newblock {A short note on the frame set of odd functions}.
\newblock {\em Bulletin of the Australian Mathematical Society}, 98(3):481–493, 2018.
\href{https://doi.org/10.1017/S0004972718000746}{10.1017/S0004972718000746}

\bibitem{Faulhuber_Hexagonal_2018}
Markus Faulhuber.
\newblock {Minimal Frame Operator Norms Via Minimal Theta Functions}.
\newblock {\em Journal of Fourier Analysis and Applications}, 24(2):545–559, 2018.
\href{https://doi.org/10.1007/s00041-017-9526-x}{10.1007/s00041-017-9526-x}

\bibitem{Fau21-RJ}
Markus Faulhuber.
\newblock An application of hypergeometric functions to heat kernels on
  rectangular and hexagonal tori and a ``{Weltkonstante}''-or-how {Ramanujan}
  split temperatures.
\newblock {\em The Ramanujan Journal}, 54(1): 1--27, 2021.
\href{https://doi.org/10.1007/s11139-019-00224-2}{10.1007/s11139-019-00224-2}

\bibitem{FauSha22}
Markus Faulhuber and Irina Shafkulovska.
\newblock{Gabor frame bound optimizations},
\newblock{\em arXiv preprint}, 2022.
\href{https://arxiv.org/abs/2204.02917}{arXiv:2204.02917}

\bibitem{FauSte17}
Markus Faulhuber and Stefan Steinerberger.
\newblock {Optimal Gabor frame bounds for separable lattices and estimates for Jacobi theta functions}.
\newblock {\em Journal of Mathematical Analysis and Applications}, 445(1):407--422, 2017.
\href{https://doi.org/10.1016/j.jmaa.2016.07.074}{10.1016/j.jmaa.2016.07.074}

\bibitem{FeiLue06}
Hans G.~Feichtinger and Franz Luef.
\newblock{Wiener amalgam spaces for the fundamental identity of Gabor analysis}.
\newblock{Collectanea Mathematica}, 57(extra):233--253, 2006.
\url{http://eudml.org/doc/41784}

\bibitem{FeiZim98}
Hans~G.\ Feichtinger and Georg Zimmermann.
\newblock {A Banach space of test functions in Gabor analysis}.
\newblock In Hans~G. Feichtinger and Thomas Strohmer, editors, {\em {Gabor
  Analysis and Algorithms: Theory and Applications}}, page 123–170.
  Birkhäuser, 1998.
  \href{https://doi.org/10.1007/978-1-4612-2016-9_4}{10.1007/978-1-4612-2016-9\_4}

\bibitem{Fin03}
Steven~R.\ {Finch}.
\newblock {\em {Mathematical constants.}}
\newblock Cambridge: Cambridge University Press, 2003.

\bibitem{Fol89}
Gerald~B.\ Folland.
\newblock {\em {Harmonic analysis in phase space}}.
\newblock Number 122 in {Annals of Mathematics Studies}. Princeton University
  Press, 1989.
  \href{https://doi.org/10.1515/9781400882427}{10.1515/9781400882427}

\bibitem{FryNaz09}
Alexander Fryntov and Fedor Nazarov.
\newblock{New estimates for the length of the Erdös–Herzog–Piranian lemniscate}.
In A.\ Alexandrov, A.\ Baranov, and S.\ Kislaykov (eds.) \newblock{\em Linear and Complex Analysis} pp.\ 49--60, American Mathematical Society Translations Series 2, vol\. 226, 2009.
\href{https://doi.org/10.1090/trans2/226}{10.1090/trans2/226}

\bibitem{Gab46}
Dennis Gabor.
\newblock {Theory of communication}.
\newblock {\em Journal of the Institution of Electrical Engineers}, 93(26):429--457, 1946.
\href{https://digital-library.theiet.org/content/journals/10.1049/ji-3-2.1946.0074}{10.1049/ji-3-2.1946.0074}

\bibitem{Gauss_WerkeIII}
Carl~F.\ Gauss.
\newblock {\em {Werke III, Analysis}}.
\newblock Königliche Gesellschaft der Wissenschaften zu Göttingen.
  Göttingen: Universitäts-Druckerei., 1866.

\bibitem{Gos11}
Maurice A.~de Gosson.
\newblock {\em {Symplectic Methods in Harmonic Analysis and in Mathematical
  Physics}}, volume~7 of {\em {Pseudo-Differential Operators. Theory and
  Applications}}.
\newblock Birkhäuser/Springer Basel AG, Basel, 2011.
\href{https://doi.org/10.1007/978-3-7643-9992-4}{10.1007/978-3-7643-9992-4}

\bibitem{Gos17}
Maurice A.\ de Gosson.
\newblock{\em The Wigner transform}.
\newblock{World Scientific}, Singapore, 2017.
\href{https://doi.org/10.1142/q0089}{10.1142/q0089}

\bibitem{GroKop19}
Karlheinz Gr{\"o}chenig and Sarah Koppensteiner.
\newblock {Gabor Frames: Characterizations and Coarse Structure}.
\newblock In A.\ Aldroubi, C.\ Cabrelli, S.\ Jaffard, and U.\ Molter (eds.), {\em New Trends in Applied Harmonic Analysis, Volume 2},
{Applied and Numerical Harmonic Analysis}, pages 93--120. Springer, 2019.
\href{https://doi.org/10.1007/978-3-030-32353-0_4}{10.1007/978-3-030-32353-0\_4}

\bibitem{Gro01}
Karlheinz Gr\"ochenig.
\newblock {\em {Foundations of Time-Frequency Analysis}}.
\newblock {Applied and Numerical Harmonic Analysis}. Birkh\"auser, Boston, MA,
  2001.
  \href{https://doi.org/10.1007/978-1-4612-0003-1}{10.1007/978-1-4612-0003-1}

\bibitem{GroHaiRom16}
Karlheinz Gr\"ochenig, Antti Haimi, and Jos\'e~L.\ Romero.
\newblock {Completeness of Gabor systems}.
\newblock {\em Journal of Approximation Theory}, 207:283–300, 2016.
\href{https://doi.org/10.1016/j.jat.2016.03.001}{10.1016/j.jat.2016.03.001}

\bibitem{HirGarBor94}
Michael Hirschhorn, Frank Garvan, and Jon Borwein.
\newblock{Cubic analogues of the Jacobian theta function $\theta(z,q)$}.
\newblock{\em Canadian Journal of Mathematics}, 45(4):673--694, 1993.
\href{https://doi.org/10.4153/CJM-1993-038-2}{10.4153/CJM-1993-038-2} 

\bibitem{How80}
Roger Howe.
\newblock {On the Role of the Heisenberg Group in Harmonic Analysis}.
\newblock {\em Bulletin of the American Mathematical Society (New Series)}, 3(2):821--843, 1980.
\href{https://doi.org/10.1090/S0273-0979-1980-14825-9}{10.1090/S0273-0979-1980-14825-9}


\bibitem{Jan95}
Augustus J.~E.~M. Janssen.
\newblock {Duality and biorthogonality for {W}eyl-{H}eisenberg frames}.
\newblock {\em Journal of Fourier Analysis and Applications}, 1(4):403–436,
  1995.
  \href{https://doi.org/10.1007/s00041-001-4017-4}{10.1007/s00041-001-4017-4}

\bibitem{Jan96}
Augustus J.~E.~M. Janssen.
\newblock {Some Weyl-Heisenberg frame bound calculations}.
\newblock {\em Indagationes Mathematicae}, 7(2):165–183, 1996.
\href{https://doi.org/10.1016/0019-3577(96)85088-9}{10.1016/0019-3577(96)85088-9}

\bibitem{Jan98}
Augustus J.~E.~M.\ Janssen.
\newblock{The duality condition for Weyl-Heisenberg frames}. In
\newblock{\em Gabor analysis and algorithms}, Applied and Numerical Harmonic Analysis, pp.\ 33--84, Birkh\"auser Boston, Boston, MA, 1998.
\href{https://doi.org/10.1007/978-1-4612-2016-9_2}{10.1007/978-1-4612-2016-9\_2}

\bibitem{Kle03}
Felix Klein.
\newblock{Gau\ss' wissenschaftliches Tagebuch, 1796--1814}.
\newblock{\em Mathematische Annalen}, 57:1--34, 1903.
\href{https://doi.org/10.1007/BF01449013}{10.1007/BF01449013}

\bibitem{Kot33}
Vladimir~A.\ Kotelnikov.
\newblock {On the Transmission Capacity of the 'Ether' and Wire in
  Electrocommunications (1933) translated by V. E. Katsnelson}.
\newblock In John~J. Benedetto and Paulo J. S.~G. Ferreira, editors,
{\em {Modern Sampling Theorey: Mathematics and Applications}}. Birkhäsuer, 2001.
\href{https://doi.org/10.1007/978-1-4612-0143-4_2}{10.1007/978-1-4612-0143-4\_2}

\bibitem{Lan29}
Edmund Landau.
\newblock {Über die Blochsche Konstante und zwei verwandte Weltkonstanten}.
\newblock {\em Mathematische Zeitschrift}, 30(1):608–634, 1929.
\href{https://doi.org/10.1007/BF01187791}{10.1007/BF01187791}

\bibitem{Low85}
Francis E.~Low.
\newblock {Complete sets of wave packets}.
\newblock {\em A passion for physics–essays in honor of Geoofrey Chew, World Scientific, Singapore}, pp.\ 17--22, 1985.
\href{https://doi.org/10.1142/9789811219207_0005}{10.1142/9789811219207\_0005}

\bibitem{Lyu92}
Yurii Lyubarskii.
\newblock {Frames in the Bargmann space of entire functions}.
\newblock In {\em {Entire and Subharmonic Functions}}, pp.\ 167--180. American Mathematical Society, Providence, RI, 1992.
\href{https://doi.org/10.1090/advsov/011}{10.1090/advsov/011}

\bibitem{Mon88}
Hugh~L.\ Montgomery.
\newblock {Minimal theta functions}.
\newblock {\em Glasgow Mathematical Journal}, 30(1):75--85, 1988.
\href{https://doi.org/10.1017/S0017089500007047}{10.1017/S0017089500007047}

\bibitem{Mum_Tata_I}
David Mumford.
\newblock {\em {Tata Lectures on Theta I}}.
\newblock {Modern Birkäuser Classics}. Birkhäuser, 2007. \\
\href{https://doi.org/10.1007/978-1-4899-2843-6}{10.1007/978-1-4899-2843-6}

\bibitem{Neumann_Quantenmechanik_1932}
John~von Neumann.
\newblock {\em {Mathematische Grundlagen der Quantenmechanik}}.
\newblock Springer-Verlag, Berlin, 2. edition, 1996.
\href{https://doi.org/10.1007/978-3-642-61409-5}{10.1007/978-3-642-61409-5}

\bibitem{Per71}
Askold Perelomov.
\newblock {On the completeness of a system of coherent states}.
\newblock {\em Theoretical and Mathematical Physics}, 6(2):156–164, 1971.
\href{https://doi.org/10.1007/BF01036577}{10.1007/BF01036577}

\bibitem{Rad43}
Hans Rademacher.
\newblock {On the Bloch-Landau Constant}.
\newblock {\em American Journal of Mathematics}, 65(3):387–390, 1943.
\href{https://doi.org/10.2307/2371963}{10.2307/2371963}

\bibitem{RonShe97}
Amos Ron and Zuowei Shen.
\newblock {Weyl-Heisenberg frames and Riesz bases in $L_2(\mathbb{R}^d)$}.
\newblock {\em Duke Mathematical Journal}, 89(2):237–282, August 1997.
\href{https://doi.org/10.1215/S0012-7094-97-08913-4}{10.1215/S0012-7094-97-08913-4}
\bibitem{Sei92_1}
Kristian Seip.
\newblock {Density theorems for sampling and interpolation in the Bargmann --
  Fock space I}.
\newblock {\em Journal für die reine und angewandte Mathematik (Crelles
  Journal)}, 429:91--106, 1992.
  \href{https://doi.org/10.1515/crll.1992.429.91}{10.1515/crll.1992.429.91}

\bibitem{SeiWal92}
Kristian Seip and Robert Wallstén.
\newblock {Density theorems for sampling and interpolation in the Bargmann --
  Fock space II}.
\newblock {\em Journal für die reine und angewandte Mathematik (Crelles
  Journal)}, 429:107--114, 1992.
\href{https://doi.org/10.1515/crll.1992.429.107}{10.1515/crll.1992.429.107}

\bibitem{Sha49}
Claude~E.\ Shannon.
\newblock Communication in the presence of noise.
\newblock {\em Proceedings of the IRE}, 37(1):10--21, 1949.
\href{https://doi.org/10.1109/JRPROC.1949.232969}{10.1109/JRPROC.1949.232969}

\bibitem{SteSha_Complex_03}
Elias Stein and Rami Shakarchi.
\newblock {\em {Complex Analysis}}.
\newblock Princeton University Press, Princeton, NJ, 2003.
\url{https://press.princeton.edu/books/hardcover/9780691113852/complex-analysis}
\bibitem{StrBea03}
Thomas Strohmer and Scott Beaver.
\newblock {Optimal OFDM design for time-frequency dispersive channels}.
\newblock {\em Communications, IEEE Transactions}, 51(7):1111–1122, July
  2003.
\href{https://doi.org/10.1109/TCOMM.2003.814200}{10.1109/TCOMM.2003.814200}

\bibitem{TolOrr95}
Richard Tolimieri and Richard~S.\ Orr.
\newblock {Poisson summation, the ambiguity function, and the theory of Weyl-Heisenberg frames.}
\newblock {\em Journal of Fourier Analysis and Applications}, 1(3):233--247, 1995.
\href{https://doi.org/10.1007/s00041-001-4011-x}{10.1007/s00041-001-4011-x}

\bibitem{WexRaz90}
Jason Wexler and Shalom Raz.
\newblock {Discrete Gabor expansions}.
\newblock {\em Signal Processing}, 21(3):207–220, November 1990.
\href{https://doi.org/10.1016/0165-1684(90)90087-F}{10.1016/0165-1684(90)90087-F}

\bibitem{Whi15}
Edmund~T.\ Whittaker.
\newblock {On the Functions which are represented by the Expansions of the Interpolation-Theory}.
\newblock {\em Proceedings of the Royal Society of Edinburgh}, 35:181–194, 1915.
  \href{https://doi.org/10.1017/S0370164600017806}{10.1017/S0370164600017806}

\bibitem{WhiWat69}
Edmund~T.\ Whittaker and George~N.\ Watson.
\newblock {\em {A Course of Modern Analysis}}.
\newblock Cambridge University Press, reprinted edition, 1969.
\href{https://doi.org/10.1017/CBO9780511608759}{10.1017/CBO9780511608759}
\bibitem{Zak67}
Joshua Zak.
\newblock {Finite Translations in Solid-State Physics}.
\newblock {\em Phys. Rev. Lett.}, 19(24):1385–1387, 1967.
\href{https://doi.org/10.1103/PhysRevLett.19.1385}{10.1103/PhysRevLett.19.1385}
\end{thebibliography}

\section*{References}
\renewcommand\refname{}

\end{document}